\RequirePackage{amsthm}
\documentclass[pdflatex,sn-mathphys-num]{sn-jnl}% Math and Physical Sciences Numbered Reference Style 
\usepackage{graphicx}%
\usepackage{multirow}%
\usepackage{amsmath,amssymb,amsfonts}%
\usepackage{mathrsfs}%
\usepackage{appendix}
\usepackage[dvipsnames]{xcolor}
\usepackage{xcolor}%
\usepackage{textcomp}%
\usepackage{manyfoot}%
\usepackage{booktabs}%
\usepackage{algorithm}%
\usepackage{algorithmicx}%
\usepackage{algpseudocode}%
\usepackage{listings}%
\usepackage[caption=false]{subfig}
\usepackage[numbers,sort&compress]{natbib}
\usepackage{cleveref}
\usepackage{soul}

\DeclareMathOperator*{\argmin}{argmin}

\definecolor{orange}{RGB}{255,127,0}

\newcommand{\remove}[1]{{\textcolor{red}{}}}
\newcommand{\crossout}[1]{{\textcolor{red}{\st{}}}}

\newcommand{\bfb}{\mathbf b}

\newcommand{\bfw}{\mathbf w}

\newcommand{\bfbtrue}{\bfb^{\mathrm{true}}}

\newcommand{\bfs}{\mathbf s}
\newcommand{\bfq}{\mathbf q}

\newcommand{\bfr}{\mathbf r}

\newcommand{\bfu}{\mathbf u}

\newcommand{\bfx}{\mathbf x}

\newcommand{\bfy}{\mathbf y}

\newcommand{\bfh}{\mathbf h}

\newcommand{\bfe}{\mathbf e}

\newcommand{\bfxtrue}{\bfx^{\mathrm{true}}}

\theoremstyle{thmstyleone}%
\newtheorem{example}{Example}%

\theoremstyle{thmstylethree}%
\begin{document}
\setstcolor{red}

\title{GKB Methods for X-Ray Computed Tomography with an Unmatched Back Projector}

% \author[1]{\fnm{Abdulmajeed} \sur{Alsubhi} }\email{ahalsubh@iu.edu.sa}
% \author[2]{\fnm{Rosemary A.} \sur{Renaut}}\email{renaut@asu.edu}
% %\equalcont{These authors contributed equally to this work.}

% \affil[1]{\orgdiv{Department of Mathematics, Faculty of Science}, \orgname{Islamic University of Madinah}, \orgaddress{\city{Madinah}, \postcode{42351}, \country{Saudi Arabia}}}

% \affil[2]{\orgdiv{School of Mathematical and Statistical Sciences}, \orgname{Arizona State University}, \orgaddress{\city{Tempe}, \postcode{85287}, \state{AZ}, \country{USA}}}

\author[]{\fnm{Abdulmajeed} \sur{Alsubhi}}\email{ahalsubhi@iu.edu.sa}
\affil[]{\orgdiv{Department of Mathematics, Faculty of Science}, \orgname{Islamic University of Madinah}, \orgaddress{\city{Madinah}, \postcode{42351}, \country{Saudi Arabia}}}

%\affil[3]{\orgdiv{Department}, \orgname{Organization}, \orgaddress{\street{Street}, \city{City}, \postcode{610101}, \state{State}, \country{Country}}}


\abstract{In large-scale X-ray Computed Tomography (CT) inverse problems, the forward and back projectors are often generated using different discretizations. This discrepancy leads to unmatched pairs of projections, resulting in inconsistent normal equations. Consequently, employing the Conjugate Gradient method does not produce a useful solution. For matched operator pairs, the Golub–Kahan bidiagonalization (GKB) method provides an efficient solution strategy. It works by projecting the original large-scale problem onto a lower-dimensional subspace, enabling the solution to be computed via a singular value decomposition of a sparse lower bidiagonal matrix. To address unmatched-pair problems in CT, we propose the AB- and BA-GKB algorithms as preconditioned forms of the GKB. These methods are straightforward to implement and allow for parameter tuning. We provide a discussion on the theoretical computational costs of our proposed algorithms in terms of floating point operations and compare with existing methods. While many Krylov methods tend to amplify noise in solutions, leading to semiconvergence, our proposed algorithms demonstrate greater resilience against this effect. We validate the effectiveness of our approach through numerical examples across various CT problems, showcasing its ability to deliver more stable solutions.}


%\pacs[MSC Classification]{65F22, 65F10, 68W40}

\maketitle

\section{Introduction}\label{sec:intro}
X-ray Computed Tomography (CT) works by placing an X-ray source on the side of an object, where the 
X-rays pass through the object and are recorded in the detector sensors. Discretization of this phenomenon leads to the forward mathematical model of the form
\begin{equation}\label{eq:inv prob}
    A\bfx\approx\bfb
\end{equation}
where $A\in \mathbb{R}^{m\times n}$ represents a discretization of the forward projector and $\bfx$ is the discretized image which needs to be reconstructed from the measured data $\bfb$. Here, we assume the measured data is contaminated by additive white Gaussian noise. In large-scale CT inverse problems, the matrix $A$ is often too large to store explicitly and direct operations are too expensive. To bypass these computational and storage obstacles, the matrix-vector products with $A$ and its transpose are computed by treating the forward projector and the back projector as functions. These functions employ different discretization schemes to construct the projectors based on approximations of the underlying physics. Thus, the back projector matrix $B \in \mathbb{R}^{n\times m} $ is not equivalent to the transpose of the forward projector $ A^\top$ and $[B, A] $ is called an unmatched projector pair \cite{article,hansen2022gmres}.

In the presence of unmatched matrix pairs, applying the Conjugate Gradient (CG) and Generalized Minimal Residual (GMRES) algorithms to the least-squares solution of the non-symmetric normal equations associated with \eqref{eq:inv prob} may fail to yield a meaningful reconstruction of $\bfx$ \cite{article}. Although hybrid regularization techniques—where regularization is incorporated into the projected problem within a Krylov subspace—are generally known to improve stability and accuracy \cite{ChungGazzolareview}, they cannot be directly applied to CT problems involving unmatched projectors. The impact of using an unmatched pairs on the convergence of multiple iterative algorithms was discussed in \cite{doi:10.1137/17M1133828}.  To address these challenges, the preconditioned GMRES, namely, the AB- and BA-GMRES methods were introduced in \cite{doi:10.1137/070696313}. These methods aim to improve the quality of the solutions of inverse problems with unmatched back projectors, by employing a preconditioner matrix $B$ constructed through appropriate diagonal scaling. Subsequently, the use of unmatched projectors in CT imaging was further investigated in \cite{hansen2022gmres}, where $B$ represents the back projector. In this framework, the standard GMRES is applied to the unmatched normal equations, in which the square systems ($AB$ and $BA$) are formed without explicitly constructing the adjoint operator. While GMRES converges rapidly, it exhibits semi-convergence: the solution converges toward the desired solution in early iterations, then diverges toward a noisy solution as iterations continue. The GMRES method computes approximations by minimizing the residual norm over a Krylov subspace generated via the Arnoldi process, leading to a reduced Hessenberg system \cite{saad2003iterative,gazzola2019arnoldi}. Consequently, the quality of the reconstruction of $\bfx$ in \eqref{eq:inv prob} is closely tied to the properties of this projected Hessenberg system in the AB- and BA-GMRES formulations.

Gene H. Golub and William Kahan introduced the GKB method in \cite{golub1965calculating} as a procedure for reducing a general matrix to bidiagonal form, which is particularly useful for computing singular value decomposition (SVD) and solving least-squares problems. The method also underpins many Krylov subspace techniques and forms the basis of modern iterative solvers such as LSQR \cite{paige1982lsqr} and LSMR \cite{fong2011lsmr}. Further applications of GKB algorithms can be explored in works such as \cite{bjorck1988bidiagonalization,bentbib2018solution,beik2020golub,bianchi2026iterated}. Motivated by these developments, we address CT problems with unmatched projectors by proposing the AB- and BA-GKB methods.  As the resulting lower bidiagonal matrix produced by GKB is sparser than the Hessenberg matrix generated by GMRES and retains important numerical properties of the original problem, the proposed approach has the potential to yield accurate reconstructions.

\subsection{Major Contributions}\label{sec:major}
We propose the AB- and BA-GKB algorithms for handling unmatched projector pairs arising in X-ray CT problems. We discuss the theoretical computational complexity of the proposed algorithms and compare them with the AB- and BA-GMRES methods. Our findings indicate that the proposed algorithms are slower than GMRES methods in terms of floating-point operations (FLOPS). Moreover, the relative cost between AB- and BA- formulations depends on whether the system is overdetermined or underdetermined. We investigate the use of two stopping criteria for iterations, including the Discrepancy Principle (DP) and Residual norm Stagnation (RNS). Through various CT problems, we demonstrate that our proposed algorithm effectively provides solutions with less severe semiconvergence than those obtained with the AB- and BA-GMRES algorithms.

\subsection{Organization }\label{sec:organ}
The organization of this paper is outlined as follows. In \Cref{sec:Background}, we review the mathematical techniques used to solve inverse problems with unmatched pairs using GMRES. In \Cref{sec:GKB}, we consider the GKB algorithm and introduce the AB- and BA-GKB. In \Cref{sec:comp costs}, we discuses the computational complexity of the proposed methods in terms of FLOPS and compare them with GMRES methods. \Cref{sec:Numerical examples} presents our numerical examples, while \Cref{sec:Conclusion} concludes our study.

\section{Background} \label{sec:Background}
Consider the least squares problem for \eqref{eq:inv prob}
\begin{equation}\label{eq:LS}
    \min_{\bfx \in \mathbb{R}^n}\|A\bfx-\bfb \|_2,
\end{equation}
where we assume that $A$ has no restrictions on its dimensions or its rank. If $m\geq n$,  the solution of \eqref{eq:LS} is equivalent to solving the normal equations
\begin{equation} \label{eq:over NE}
   A^\top A\bfx =A^\top \bfb.
\end{equation}
When $m < n$, the solution takes the form
\begin{equation} \label{eq:under NE}
   A A^\top \bfy = \bfb, \ \ \ \bfx=A^\top \bfy.
\end{equation}
When the projectors are not matched, e.g. as in CT,   the normal equations in \cref{eq:over NE,eq:under NE} are replaced with the so-called unmatched normal equations, respectively, 
\begin{equation} \label{eq:over UN NE}
   B A\bfx =B \bfb,
\end{equation}
and
\begin{equation} \label{eq:under UN NE}
   A B \bfy = \bfb, \ \ \ \bfx=B \bfy.
\end{equation}
Noting that when $B=A^\top$, \eqref{eq:over UN NE} and \eqref{eq:under UN NE} are equivalent to \eqref{eq:over NE} and \eqref{eq:under NE}, respectively. This situation is called the matched normal equations \cite{hansen2010discrete}. If the range of $A^\top$ is equivalent to the range of $B$, then the matched and unmatched normal equations are equivalent \cite{doi:10.1137/070696313}. The solutions in this case can be computed using CGLS or LSQR. However, our work specifically addresses the scenario involving unmatched projectors, where $B \neq A^\top$.

In the following sections, we review AB- and BA-GMRES in relation to previous studies on unmatched projectors.

\subsection{AB-GMRES and BA-GMRES} \label{sec:ABBA GMRES}
GMRES is an efficient iterative Krylov subspace method for solving a square linear least-squares problem \cite{saad2003iterative}. This algorithm terminates when the $\ell_2$ norm of the residual becomes sufficiently small, or falls below a specified threshold. In many CT problems, the system matrix $A$ is not square, so applying GMRES directly to $A$ is not feasible. It is important to note that the system matrices $AB$ and $BA$ in \eqref{eq:over UN NE} and \eqref{eq:under UN NE} are square but generally non-symmetric, so the popular CGLS is not applicable in such a case. The AB-GMRES and BA-GMRES algorithms, introduced in \cite{doi:10.1137/070696313}, address this dimensionality issue by allowing GMRES to be applied to rectangular systems, employing the square matrices $AB$ and $BA$ implicitly without the need for explicit multiplications.

AB-GMRES applies GMRES to \eqref{eq:under UN NE}, in which B serves as the right preconditioner. At iteration $k$, the AB-GMRES computes the solution $\bfx_k=B\bfy_k$ that minimizes $\|AB\bfy-\bfb \|^2_2$. Let $\bfx_0$ be an initial guess of the solution $\bfx$, then the initial residual is computed by $\bfr_0=\bfb-A\bfx_0$. A related approach was presented in \cite{calvetti2000gmres} to solve an overdetermined system $(m\geq n)$ using GMRES. In particular, the authors extend the columns of $A$ on the right side by incorporating zero padding, resulting in a new square system matrix 
\begin{equation*}
\tilde{A}=
    \begin{bmatrix}
        A & 0_{m\times (m-n)}
    \end{bmatrix}
    \in \mathbb{R}^{m\times m}
\end{equation*} and then applied the GMRES to the new least squares problem 
\begin{equation*}
    \min_{\bfy\in \mathbb{R}^m} \|\tilde{A} \bfy -\bfb\|_2 = \min_{\bfx\in \mathbb{R}^n} \|A \bfx -\bfb\|_2.
\end{equation*}
This approach can be viewed as a particular case of AB-GMRES, where $\tilde{A} =AB$ and $B=[I_n,0]\in \mathbb{R}^{n\times m}$.

BA-GMRES applies GMRES to \eqref{eq:over UN NE}, with B serving as the left preconditioner. At iteration $k$, the BA-GMRES computes the solution $\bfx$ that minimize $\|BA\bfx-B\bfb \|^2_2$. The initial residual in this case is computed by $\bfr_0=B( \bfb-A\bfx_0)$. A related approach was proposed in \cite{reichel2005breakdown} for an underdetermined $m < n$ system using GMRES. In particular, the authors extend the rows of $A$ from below by incorporating zero padding, resulting in a new square system 
\begin{equation}
\tilde{A}=
    \begin{bmatrix}
        A\\0
    \end{bmatrix}
    \in \mathbb{R}^{n\times n} \, \ \ \ \ \ \tilde{\bfb}=\begin{bmatrix}
    \bfb\\0
\end{bmatrix} \in \mathbb{R}^{n}.
\end{equation}  
Then, the GMRES is applied to the new least squares problem 
\begin{equation*}
    \min_\bfx \|\tilde{A} \bfx -\tilde{\bfb}\|_2.
\end{equation*}
In fact, this is a special case of BA-GMRES with 
\begin{equation*}
    B=\begin{bmatrix}
        I_m\\0
    \end{bmatrix}\in \mathbb{R}^{n\times m},
\end{equation*}
in which $\tilde{A}=BA$ and $\tilde{\bfb}=B\bfb$.

It is worth noting that both AB-GMRES and BA-GMRES use the same Krylov space $\mathcal{K}_k(BA, B\bfb)$ \cite{doi:10.1137/070696313}. Arnoldi iteration is used within the GMRES to construct the orthogonal basis needed to compute an approximate solution with minimal residual. Let $M$ be the coefficient matrix in $AB$ or $BA$, then $k$ applications of Arnoldi iteration yield
\begin{equation}\label{eq:Arn iter}
    M \bar{Q}_k=\bar{Q}_{k+1}H_k,
\end{equation}
where $\bar{Q}_{k+1}= [\bar{\bfq}_1,\bar{\bfq}_2,\dots,\bar{\bfq}_{k+1}]\ $ and the matrix $H_{k+1,k}$ is an upper Hessenberg matrix,  defined by
\begin{equation}\label{upper Hessenberg}
     H_k=
\begin{bmatrix}
   h_{1,1}&h_{1,2} &\dots & h_{1,k}\\
   h_{2,1}&h_{2,2} &\dots & h_{2,k}\\
   &\ddots &\ddots &\vdots\\
   &&h_{k,k-1}&h_{k,k}\\
   &&&h_{k+1,k}
\end{bmatrix}.
\end{equation}
Applying the decomposition \eqref{eq:Arn iter} to the minimization problems \eqref{eq:over UN NE} and \eqref{eq:under UN NE} yields the small projected problem 
\begin{equation}\label{eq:proj}
    \bfu_k=\argmin_{\bfu} \Big\| H_k \bfu - \bfe_1 \|\bfr_0\|_2 \Big\|_2,
\end{equation}
where $\bfe_1 =(1,0,\dots,0)^\top \in \mathbb{R}^{k+1}$. This small projected problem can be solved iteratively using the CGLS method. Alternatively, the solution $\bfu_k$ can be derived from matrix decompositions, such as QR, LU, and SVD. 

The steps of the AB-GMRES and BA-GMRES Algorithms are summarized in \Cref{alg:ABGMRES,alg:BAGMRES}.

\begin{minipage}{0.47\textwidth}
\begin{algorithm}[H]
    \centering
    \caption{AB-GMRES \cite{doi:10.1137/070696313}}\label{alg:ABGMRES}
    \begin{algorithmic}[1]
       \Require $A$, $B$, $\bfb$, 
       \State \text{Choose initial $\bfy_0$}, \text{set ${\bfx}_0=B\bfy_0$}
        \State {$\bfr_0=\bfb-A\bfx_0$}
        \State {$\bar{\bfq}_1=\bfr_0 / \| \bfr_0\|_2$}
        \For{$k=1,2,\dots$}
        \State ${\bar{\bfw}_k=A(B \bar{\bfq}_k})$ 
        \For{$i=1,2,\dots,k$}
        \State ${\bfh_{i,k} = \bar{\bfq}^\top_i \bar{\bfw}_k}$ 
        \State ${\bar{\bfw}_k= \bar{\bfw}_k- \bfh_{i,k} \bar{\bfq}_i }$ 
         \EndFor
        \State ${\bfh_{k+1,k} = \|\bar{\bfw}_k \|_2}$ 
        \State ${\bar{\bfq}_{k+1}=\bar{\bfw}_k / \bfh_{k+1,k}}$ 
        \State {update $\bfu_k$: solve \eqref{eq:proj}}
        \State {$\bfx_k=\bfx_0+B \bar{Q}_k \bfu_k$}
        \State \text{$\bfr_k=\bfb-A\bfx_k$}
        \State \text{Stopping rule is applied here}
        \EndFor
    \end{algorithmic}
\end{algorithm}
\end{minipage}
\hfill
\begin{minipage}{0.47\textwidth}
\begin{algorithm}[H]
    \centering
    \caption{BA-GMRES \cite{doi:10.1137/070696313}}\label{alg:BAGMRES}
    \begin{algorithmic}[1]
      \Require $A$; $B$, $\bfb$, 
        \State \text{Choose initial $\bfx_0$} 
        \State {$\bfr_0=B(\bfb-A\bfx_0)$}
        \State {$\bar{\bfq}_1=\bfr_0 / \| \bfr_0\|_2$}
        \For{$k=1,2,\dots$}
        \State ${\bar{\bfw}_k=B (A \bar{\bfq}_k})$ 
        \For{$i=1,2,\dots,k$}
        \State ${\bfh_{i,k} = \bar{\bfq}^\top_i \bar{\bfw}_k}$ 
        \State ${\bar{\bfw}_k= \bar{\bfw}_k- \bfh_{i,k} \bar{\bfq}_i }$ 
         \EndFor
        \State ${\bfh_{k+1,k} = \|\bar{\bfw}_k \|_2}$ 
        \State ${\bar{\bfq}_{k+1}=\bar{\bfw}_k / \bfh_{k+1,k}}$ 
        \State {update $\bfu_k$: solve \eqref{eq:proj}}
        \State {$\bfx_k=\bfx_0+ \bar{Q}_k \bfu_k$}
        \State \text{$\bfr_k=B(\bfb-A\bfx_k)$}
        \State \text{Stopping rule is applied here}
        \EndFor
    \end{algorithmic}
\end{algorithm}
\end{minipage}
% The approximate solution associated with \cref{eq:proj} for AB-GMRES is given by
% \begin{equation}\label{eq:approx pro AB}
%     \bfx_k=\bfx_0+B \bar{Q}_k \bfu_k,
% \end{equation}
% as well as for BA-GMRES
% \begin{equation}\label{eq:approx pro BA}
%     \bfx_k=\bfx_0+ \bar{Q}_k \bfu_k.
% \end{equation}

\vspace{0.5cm} % I use this manually

It is worth emphasizing that the AB-GMRES algorithm explicitly requires the application of the operator $B \bar{Q}_k$ at each iteration. In contrast, the BA-GMRES algorithm incorporates this operation implicitly through the residual computation, thereby avoiding its direct evaluation. This difference affects not only the computational complexity of these algorithms but also their convergence. A significant drawback of the GMRES method is that the growth of the Krylov subspace dimension, thereby increasing storage demands \cite{saad2003iterative}. To address this challenge and enhance the practicality of the method, one can implement GMRES with restarting to reduce the total number of iterations needed for convergence. However, it is important to highlight that restarted GMRES can stagnate and exhibit slower convergence, particularly if the matrix is not positive definite \cite{saad2003iterative}. Therefore, we use the full GMRES, which guarantees convergence.

The AB- and BA-GMRES methods inherit the convergence behavior of Krylov subspace methods and exhibit semiconvergence when applied to ill-conditioned inverse problems. Specifically, during the initial iterations, the computed solution approaches the true solution as the Krylov subspaces capture the dominant singular components associated with the largest singular values. However, with further iterations, the subspaces begin to include components corresponding to smaller singular values, which are typically corrupted by noise. This semiconvergence issue can worsen with unmatched projectors because the combined operators $AB$ and $BA$ tend to be non-normal and exhibit spectral properties different from those observed in matched normal equations. To tackle this challenge, stopping rules can be applied to terminate the iteration when the convergence behavior changes. The stopping rules in AB- and BA-GMRES, provided in \Cref{alg:ABGMRES,alg:BAGMRES}, serve as an implicit regularization \cite{bjorck1996numerical}.

The SVD is one of the most powerful and reliable techniques for solving inverse problems. Let  $H_k=U\Sigma V^\top$ denote the SVD of the upper Hessenberg matrix in \eqref{upper Hessenberg} arising from the projected system \eqref{eq:proj}, where $U$ and $V$ are two unitary matrices and $\Sigma$ is a diagonal matrix containing the singular values of $H_k$ in decreasing order. These SVD components facilitate the computation of the iterative solution  $\bfu$ at each iteration. It is crucial to note that the quality of the computed solution $\bfx$ in both AB- and BA-GMRES methods is significantly influenced by the condition of the relatively sparse matrix $H_k$. From a computational standpoint, the sparsity structure of $H_k$ directly affects the efficiency and applicability of GMRES. In particular, the upper triangular part in $H_k$, including the main diagonal, has $k(k-1)/2$ entries, while the lower subdiagonal has $k-1$ entries. Thus, the total number of possible nonzero entries in $H_k$ is roughly $(k^2+3k)/2$ elements. As the size of the projected problem increases, computing the SVD of $H_k$ might be computationally expensive. Therefore, in the following section, we consider an alternative approach based on a well-established decomposition that yields a projected system matrix with greater sparsity than $H_k$.

\section{AB- and BA-GKB} \label{sec:GKB}
In this section, we introduce AB-GKB and BA-GKB algorithms that incorporate unmatched projectors into the GKB process. To comprehend our proposed methods better, we first review the GKB algorithm. 

\subsection{GKB with a Matched Transpose}
The GKB algorithm is an iterative technique introduced in \cite{golub1965calculating} and used to reduce a large matrix to a lower bidiagonal form, with the aim of solving ill-posed inverse problems. After $k$ steps of the GKB with initial vector $\bfb$, the matrix $A$ can be factorized into 
\begin{equation}\label{eq:GKB dec}
    A Q_k=S_{k+1} C_k,
\end{equation}
where $Q_{k}$ and $S_{k+1}$ are column orthogonal matrices and
\begin{equation}\label{eq:lower bid}
   C_{k} = \begin{bmatrix}
       \alpha_1&&&&\\
   t_2&\alpha_2&&&\\
   &\ddots&\ddots&&\\
      &&&t_k&\alpha_k\\
     &&&&t_{k+1}
   \end{bmatrix}
\end{equation} 
is a lower bidiagonal matrix. The GKB generates the two Krylov subspaces $\mathcal{K}_k(A^\top A, A^\top \bfb)$ and $\mathcal{K}_k(A A^\top, \bfb)$. 
This approach enables the reduction of the inverse problem in \eqref{eq:inv prob} to a significantly smaller lower bidiagonal system. The approximate solution within this projected subspace can then be computed efficiently using the SVD \cite[Ch. 7]{hansen1998rank}. From a spectral viewpoint, the singular spectrum of the $C_k$ in \eqref{eq:lower bid} includes elements that are good approximations to the biggest dominant singular values of $A$ \cite{paige1982lsqr}. Consequently, computing the SVD of $C_k$ yields an approximate truncated SVD (TSVD) of the original matrix $A$. The steps of the GKB method for approximating the solution to the inverse problem \eqref{eq:inv prob} are summarized in \Cref{alg:GKB}.

\begin{algorithm}[H] 
    \centering
    \caption{GKB}\label{alg:GKB}
    \begin{algorithmic}[1]
 \Require $A$,  $\bfb$
       \State $t_1=\|\bfb\|_2$ 
       \State $\bfs_1=\bfb/t_1$, \ $\bfq=A^\top \bfs_1$
       \State $\alpha_1=\|\bfq\|_2$
       \State $\bfq_1=\bfq/\alpha_1$
        \For{$k=1,2,\dots$}
\State {$\bfs_k=A \bfq_{k-1}- \alpha_{k-1} \bfs_{k-1}$ }
\State {$t_k=\|\bfs_k \|_2$, \ $\bfs_k=\bfs_k /t_k $ }
\State {$\bfq_k= A^\top \bfs_k- t_k \bfq_{k-1}$ \;}
\State {$\bfq_k = \bfq_k - \sum^{k-1}_{i=1} (\bfq^\top_i \bfq_k)  \bfq_i$}
\State {$\alpha_k=\|\bfq_k \|_2, \ \bfq_k=\bfq_k/\alpha_k$ \;}
\State{$\bfs_i$ and $\bfq_i$ define $S_{k+1}$,  $Q_k$}
\State{Scalars $t_k$ and $\alpha_k$ define  $C_k$}       
\State \text{Find $\bfu_k$ which minimizes} ${\|C_k \bfu - t_1 \bfe^{k+1}_1\|_2}$
        \State {$\bfx_k= Q_k \bfu_k$}
        \State \text{Stopping rule is applied here}
        \EndFor
    \end{algorithmic}
\end{algorithm}

  One-sided reorthogonalization is applied at step $9$ of \Cref{alg:GKB} to maintain the orthogonality of the columns of $Q$ \cite{simon2000low}. It is evident that the implementation of the GKB is equivalent to LSQR, which proceeds by applying the CG to the normal equations \eqref{eq:over NE} \cite{paige1982lsqr}. To improve the TSVD approximation of AAA and better capture the tail of the spectrum, one may compute a higher-rank approximation and subsequently truncate it to rank $k$\cite{larsen1998lanczos,renaut2017hybrid,alsubhi2025enlarged}. However, determining an appropriate size for the enlarged Krylov subspace is nontrivial: underestimation may degrade the reconstruction quality, whereas overestimation can significantly increase computational cost. Therefore, we avoid enlarging the Krylov subspace and instead terminate the GKB iterations using an appropriate stopping criterion.

\Cref{alg:GKB} is more general than full GMRES and can be applied to any system matrix without restrictions on its dimensions. Each iteration of \Cref{alg:GKB} involves evaluating both matrix-vector multiplication with $A$ and its transpose in a matrix-free fashion. However, in the case of an unmatched projector pair, \Cref{alg:GKB} cannot be applied directly. To address this, the next section considers a modification of the GKB process in which the matrix $B$ is used in place of $A^\top$.

\subsection{AB- and BA-GKB  with Mismatched Operators}
In this section, we present the AB- and BA-GKB algorithms for CT reconstruction with unmatched operator pairs. These methods apply GKB to the unmatched normal equations \cref{eq:over UN NE,eq:under UN NE}, while avoiding the explicit formation of the composite operators$AB$ and $BA$.

AB-GKB applies GKB to \eqref{eq:under UN NE}, in which $B$ serves as the right preconditioner. At iteration $k$, the AB-GKB reduces the large scale problems
\begin{equation*}
    \min_{\bfx\in \mathbb{R}^n}\|A B \bfx-\bfb\|_2^2
\end{equation*} 
to the following small projected problem
\begin{equation} \label{eq:proj GKB}
\min_{\bfu \in \mathbb{R}^k}\|C_k\bfu- t_1 \bfe^{k+1}_1\|_2^2.    
\end{equation}
The approximate solution is then recovered by projection onto the original space, yielding $\bfx = B Q_k \bfu_k$. In this formulation, the starting vector is analogous to that used in the standard GKB method described in \Cref{alg:GKB}.

In contrast, BA-GKB applies the GKB process to \eqref{eq:over UN NE}, with $B$ serving as a left preconditioner. At iteration $k$, BA-GKB transforms the problem
\begin{equation*}
    \min_{\bfx\in \mathbb{R}^n}\|BA \bfx-B \bfb\|_2^2
\end{equation*}
into a projected problem of the same form as \eqref{eq:proj GKB}. The corresponding approximate solution is obtained as $\bfx = Q_k \bfu_k$. Here, the initial vector is computed as $t_1 = \| B \bfb \|_2$.

Steps of AB-GKB and BA-GKB are summarized in \Cref{alg: AB-GKB,alg: BA-GKB}.
\begin{minipage}{0.47\textwidth}
\begin{algorithm}[H]
    \centering
    \caption{AB-GKB}\label{alg: AB-GKB}
    \begin{algorithmic}[1]
 \Require $A$, $B$,  $\bfb$
       \State $t_1=\|\bfb\|_2$ 
       \State $\bfs_1=\bfb/t_1$, \ $\bfq=B^\top (A^\top \bfs_1)$
       \State $\alpha_1=\|\bfq\|_2$
       \State $\bfq_1=\bfq/\alpha_1$
        \For{$k=1,2,\dots$}
\State {$\bfs_k=A(B\bfq_{k-1})- \alpha_{k-1} \bfs_{k-1}$ }
\State {$t_k=\|\bfs_k \|_2$, \ $\bfs_k=\bfs_k /t_k $ }
\State {$\bfq_k=B^\top (A^\top \bfs_k)- t_k \bfq_{k-1}$ \;}
\State {$\bfq_k = \bfq_k - \sum^{k-1}_{i=1} (\bfq^\top_i \bfq_k)  \bfq_i$}
\State {$\alpha_k=\|\bfq_k \|_2, \ \bfq_k=\bfq_k/\alpha_k$ \;}
\State{$\bfs_i$ and $\bfq_i$ define $S_{k+1}$,  $Q_k$}
\State{Scalars $t_k$ and $\alpha_k$ define  $C_k$}       
 \State {update $\bfu_k$: solve \eqref{eq:proj GKB}}
\State {$\bfx_k=B Q_k \bfu_k$}
\State \text{Stopping rule is applied here}
\EndFor
\end{algorithmic}
\end{algorithm}
\end{minipage}
\hfill
\begin{minipage}{0.47\textwidth}
\begin{algorithm}[H]
    \centering
    \caption{BA-GKB}\label{alg: BA-GKB}
    \begin{algorithmic}[1]
                \Require $A$, $B$,  $\bfb$
       \State $t_1=\|B\bfb\|_2$ 
       \State $\bfs_1=B\bfb/t_1$, \ $\bfq=A^\top (B^\top  \bfs_1)$
       \State $\alpha_1=\|\bfq\|_2$
       \State $\bfq_1=\bfq/\alpha_1$
        \For{$k=1,2,\dots$}
\State {$\bfs_k=B(A\bfq_{k-1})- \alpha_{k-1} \bfs_{k-1}$ }
\State {$t_k=\|\bfs_k \|_2$, \ $\bfs_k=\bfs_k /t_k $ }
\State {$\bfq_k=A^\top (B^\top \bfs_k)- t_k \bfq_{k-1}$ \;}
\State {$\bfq_k = \bfq_k - \sum^{k-1}_{i=1} (\bfq^\top_i \bfq_k)  \bfq_i$}
\State {$\alpha_k=\|\bfq_k \|_2, \ \bfq_k=\bfq_k/\alpha_k$ \;}
\State{$\bfs_i$ and $\bfq_i$ define $S_{k+1}$,  $Q_k$}
\State{Scalars $t_k$ and $\alpha_k$ define  $C_k$}       
 \State {update $\bfu_k$: solve \eqref{eq:proj GKB}}
\State {$\bfx_k=Q_k \bfu_k$}
\State \text{Stopping rule is applied here}
\EndFor
    \end{algorithmic}
\end{algorithm}
\end{minipage}

\vspace{0.5cm} % I use this manually

In steps $13$ of \Cref{alg: AB-GKB,alg: BA-GKB}, the solutions of the small least squares problems can be computed efficiently using the SVD of the small lower bidiagonal matrix $C_k$. The efficiency of the AB- and BA-GKB depends heavily on the properties of $C_k$, whose dimension increases with the iteration count.  It is important to note that AB-GKB requires the multiplication $BQ_k$ at each iteration, whereas BA-GKB performs this multiplication implicitly through its residual. This difference significantly affects the computational complexity of these algorithms and their convergence properties. As with many Krylov subspace methods, GKB may exhibit semiconvergence for ill-posed inverse problems. This issue may be even more pronounced with AB-GKB and BA-GKB, as the operators $AB$ and $BA$ are non-normal and exhibit spectral properties that differ from those of the corresponding normal equations. Therefore, we consider different stopping criteria—specifically, the DP and SNR methods—to halt the iterations before the noise adversely affects the solution.

\section{Complexity of AB- and BA Algorithm} \label{sec:comp costs}
In this section, we present the computational complexity in terms of FLOPS for reconstructing CT problems using the AB-GMRES, BA-GMRES, AB-GKB, and  BA-GKB methods. We assume that operations $AB \bfy$ and $BA\bfy$ and their transposes are implemented by performing two successive matrix-vector multiplications rather than matrix-matrix multiplications.

We first comment on the computational costs of the AB-GMRES and BA-GMRES, presented in \Cref{alg:ABGMRES,alg:BAGMRES}. At iteration $k$ in both algorithms, step $5$ uses two matrix-vector multiplications, requiring $\mathcal{O}(4mnk)$ FLOPS. Steps $6$ to $9$ apply modified Gram-Schmidt (MGS) to reorthogonalize the column $\bar{\bfw}_k$ against all previous columns, which can be rewritten as $\bar{\bfw}_k = \bar{\bfw}_k - \sum^k_{i=1} \langle \bar{\bfw}_k , \bar{\bfq}_i\rangle \bar{\bfq}_i$. Since the vectors $\bar{\bfw}_k$ and $\bar{\bfq}_k$ have different sizes in the AB- and BA-GMRES algorithms, the MGS cost is approximately $\mathcal{O}(\sum^k_{i=1} 4mi) \approx 2 m k^2$ FLOPS for AB-GMRES and $\mathcal{O}(\sum^k_{i=1} 4ni) \approx 2 n k^2$ FLOPS for BA-GMRES. At step $13$ of \Cref{alg:ABGMRES}, the product $B\bar{Q}_k \bfu_k$ is computed explicitly at each iteration, requiring $\mathcal{O}(2mnk+2mk^2)$ flops. In contrast, \Cref{alg:BAGMRES} evaluates this product implicitly through the residual, reducing the cost to $\mathcal{O}(2nk^2)$ FLOPS.

For the projected problem, computing the Lanczos SVD of a sparse matrix requires $\mathcal{O}(T k nnz)$ FLOPS \cite{larsen1998lanczos},  where $T$ denotes the truncated parameter of the singular values kept at iteration $k$ and $nnz$ denotes the number of nonzero elements in the sparse matrix. Since no truncation is applied here, the cost becomes $\mathcal{O}(k^2 nnz)$. As discussed in \Cref{sec:ABBA GMRES}, the matrix $H_k$ contains approximately $(k^2+3k)/2$ entries, and its SVD requires $\mathcal{O}( (k^4+3k^2)/2)$ FLOPS. Lower-order operations are neglected.
Combining these steps together, the dominant computational complexity of the AB-GMRES method for  $k (\ll min(m,n))$ is given by
\begin{equation}\label{cost:ABGMRES}
   \mathcal{O}(6mnk+4mk^2+(k^4+3k^2)/2).
\end{equation}
Similarly, the dominant cost of the BA-GMRES algorithm is 
\begin{equation}\label{cost:BAGMRES}
    \mathcal{O}(4mnk+4nk^2 + (k^4+3k^2)/2 ).
\end{equation}

We next consider the computational cost of the AB- and BA-GKB algorithms in \Cref{alg: AB-GKB,alg: BA-GKB}. The dominant cost in both methods arises from the successive applications of $AB$, $BA$, and their transposes in steps $6$ and $8$, totaling $\mathcal{O}(8mnk)$ FLOPS. Maintaining column orthogonality in step $9$ requires approximately $\mathcal{O}(2 m k^2)$  FLOPS for AB-GKB and $\mathcal{O}(2 n k^2)$  FLOPS for BA-GKB.

At step $15$ of \Cref{alg: AB-GKB}, the product $BQ_k \bfu_k$ is evaluated explicitly at each iteration, incurring $\mathcal{O}(2mnk+2mk^2)$ FLOPS, whereas the corresponding operation in \Cref{alg: BA-GKB} is computed implicitly via the residual, requiring only $\mathcal{O}(2nk^2)$ FLOPS. Since the bidiagonal matrix $C_k$ contains at most $2k$ nonzero entries, its Lanczos SVD can be computed in $\mathcal{O}(2k^3)$ FLOPS \cite{larsen1998lanczos}. Lower-order costs are neglected. Therefore, the sum of the provided costs yields the total cost of \Cref{alg: AB-GKB}
\begin{equation}\label{cost:ABGKB}
    \mathcal{O}(10mnk+4mk^2+2k^3),
\end{equation}
as well as the total cost of \Cref{alg: BA-GKB}
\begin{equation}\label{cost:BAGKB}
    \mathcal{O}(8mnk+4nk^2+ 2k^3).
\end{equation}

A comparison of the Lanczos SVD costs for $C_k$ and $H_k$ shows that the cost associated with $C_k$ is approximately $1/k$ of that of $H_k$. Moreover, the enhanced sparsity of $C_k$ in AB- and BA-GKB can yield accurate approximations of the dominant singular values, thereby improving the quality of the reconstructed solution.

If any of the methods terminate early with a small $k$, the last terms in \cref{cost:ABGMRES,cost:BAGMRES,cost:ABGKB,cost:BAGKB} become negligible. Comparing \eqref{cost:ABGMRES} with \eqref{cost:ABGKB} and \eqref{cost:BAGMRES} with \eqref{cost:BAGKB} indicates that AB-GMRES requires approximately $3/5$ of the computational cost of AB-GKB, while BA-GMRES incurs about half the cost of AB-GKB. Furthermore, when $m > n$, the AB-based methods are typically more expensive than their BA counterparts, whereas for $m \ll n$, the opposite trend holds. Overall, for a fixed number of iterations $k$, GMRES-based methods are computationally more efficient in terms of FLOPS.

\section{Numerical Experiments}\label{sec:Numerical examples}
In this section, we assess the performance of the proposed AB-GKB and BA-GKB algorithms for reconstructing CT problems with unmatched pairs of forward and back projectors. We compare the reconstructions produced by these algorithms with those derived from AB-GMRES and BA-GMRES. For all simulations, we assume the availability of both the exact image $\bfxtrue$ and the noise-free data $\bfbtrue$. To evaluate the impact of noise, we introduce additive white Gaussian noise to $\bfbtrue$. The Signal-to-Noise Ratio (SNR) is employed to quantify the noise present in $\bfb$.

To evaluate the quality of the reconstruction, we record the relative reconstruction error (RRE) between the exact image $\bfxtrue$ and the reconstructed image $\bfx^{(k)}$ at iteration $k$, as given by   
\begin{equation}\label{eq:RE}
    \text{RRE}(\bfx^{(k)})=\frac{\|\bfx^{(k)}- \bfxtrue\|_2} {\|\bfxtrue\|_2},
\end{equation}
In most practical problems, the exact image is not available, and hence using the RRE to stop the iterations at the point of semiconvergence is not possible. Here, we assume an accurate estimation of the noise level $\|\bfe\|_2=\|\bfbtrue-\bfb\|_2$ is available and we terminate the iteration using the DP \cite{morozov2012methods} as 
\begin{equation}\label{eq:DP}
    k_{DP}=\min_k{\|A\bfx_k-\bfb\|_2\leq \tau \|\bfe\|_2},
\end{equation}
where $\tau\geq 1$ is a safety factor, which is set to be $1$.

The Residual Norm Stagnation (RNS) occurs when solver stops making meaningful progress. The iteration is terminated using RNS when 
\begin{equation}\label{eq:RNS}
    k_{RNS}= \frac{ \bigr| \|A\bfx_{k-1}-\bfb\|_2 - \|A\bfx_k-\bfb\|_2 \bigr|}  {\|A\bfx_{k-1}-\bfb\|_2} < \epsilon,
\end{equation}
where the tolerance $\epsilon$ is user-defined to ensure convergence, and its appropriate choice may vary across algorithms. To facilitate a consistent comparison between the proposed methods and GMRES-based approaches, $\epsilon$ is selected empirically for each simulation—based on numerical experiments (not reported)—to minimize the RRE for both AB-GMRES and BA-GMRES.

For each problem, we compare the results obtained using the DP in \eqref{eq:DP} with those from the  RNS in \eqref{eq:RNS}. In our implementations, we record the RRE and the corresponding iteration counts. Furthermore, we report average timings across $5$ runs to validate the theoretical computational discussion presented in \Cref{sec:comp costs}. It is important to note that the costs provided in \Cref{sec:comp costs} correspond to equivalent values of $k$ across all algorithms. However, in practice, each algorithm may terminate with a different value of $k$. To ensure a fair comparison, we therefore record the runtime for all algorithms using a fixed number of $150$ iterations.

All numerical results presented were computed using MATLAB Version 2026a \cite{MATLAB} on a desktop equipped with an Intel(R) Core(TM) i5-14400F CPU, 16 GB of memory, and an NVIDIA GeForce RTX 5060 GPU. The software will be available on request to the author.

\begin{figure}[ht!]
  \centering
   \subfloat[True image]{\label{Fig:True mri}\includegraphics[width=.25\textwidth]{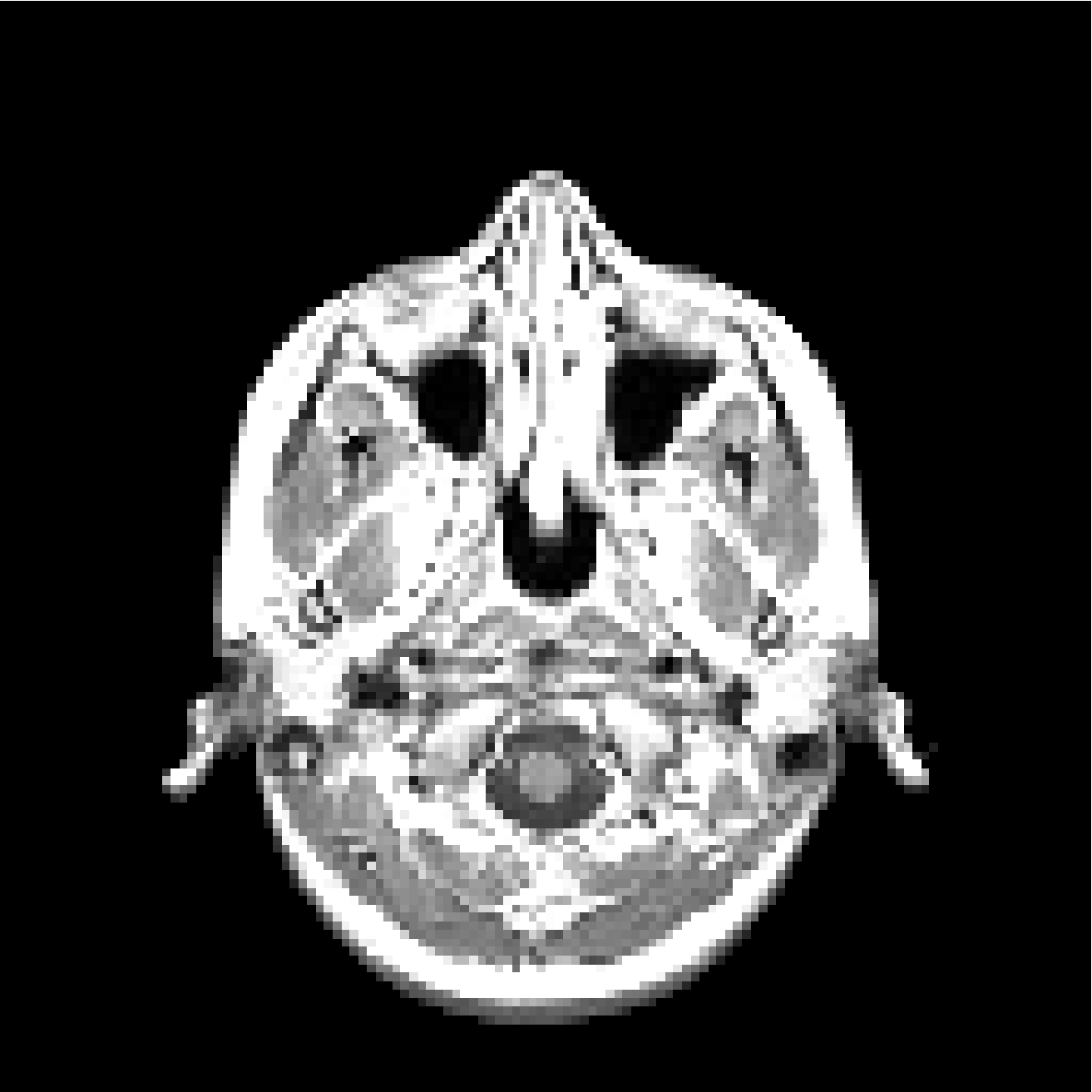}} \
   \subfloat[Noisy sinogram]{\label{Fig:blurred mri}\includegraphics[width=.25\textwidth]{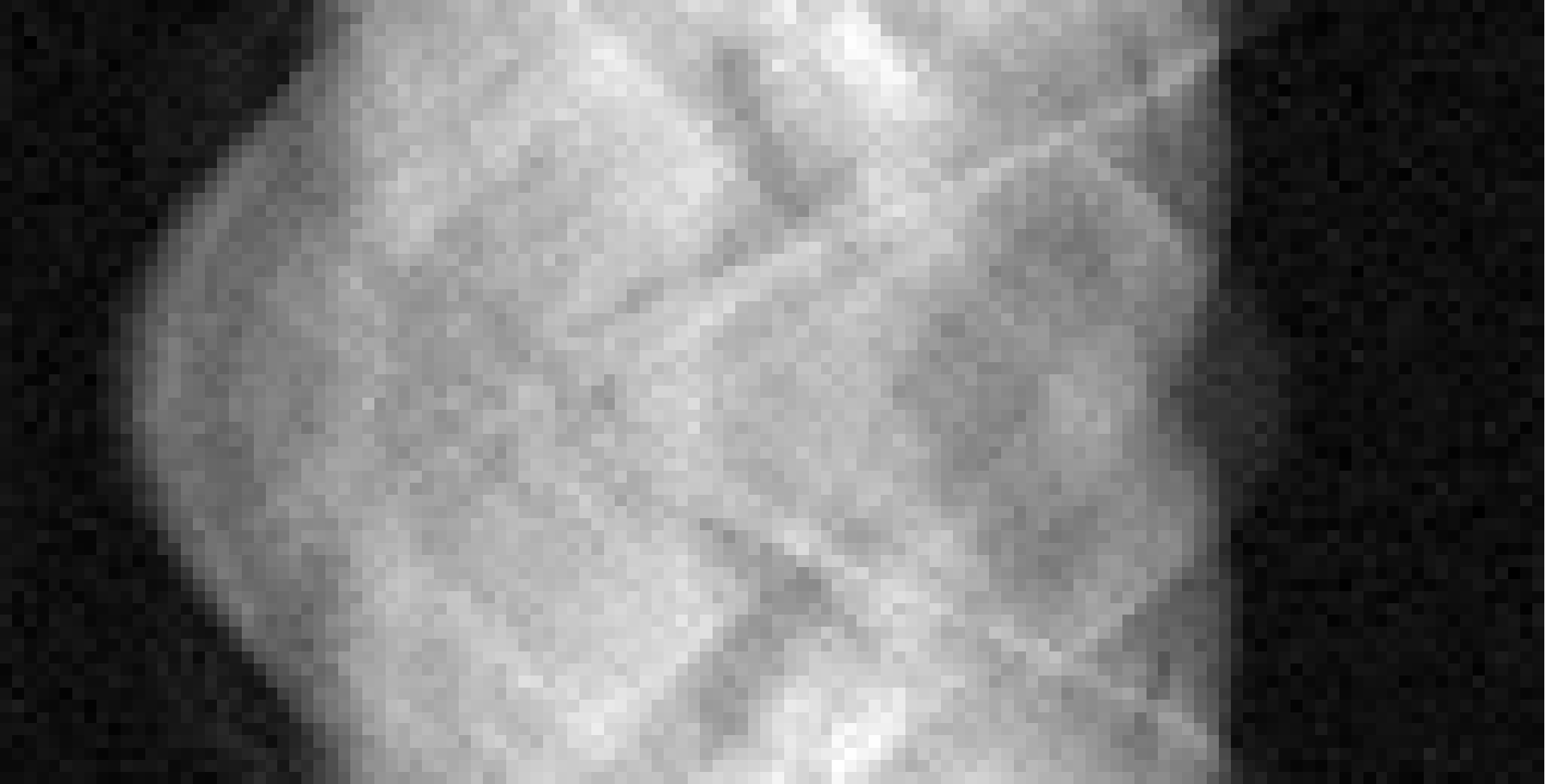}}
 \caption{The true image of $128 \times 128$ pixels and noisy sinogram with $4\%$ Gaussian noise, for the problem in \Cref{ex:problem1}. \label{Fig:mri True and noisy}}  
  \end{figure}

\begin{example}\label{ex:problem1}
We consider a phantom image of size $128 \times 128$, shown in \Cref{Fig:True mri}. The sinogram data are generated using parallel geometry through $65$ view angles. We use the ASTRA Toolbox \cite{van2015astra} to generate the forward projector $A \in \mathbb{R}^{8320\times 16384}$  and backprojector $B \in \mathbb{R}^{16384\times 8320}$.
Then we use $4 \%$ white Gaussian noise, corresponding to SNR of approximately $28$, to obtain the contaminated sinogram data $\bfb \in \mathbb{R}^{8320}$, shown in \Cref{Fig:blurred mri}.
\end{example}

This small problem is designed to test the theoretical properties of the proposed algorithms with an underdetermined coefficient matrix. When using GPU acceleration in ASTRA, the operators $A$ and $B$ are implemented as a matrix-free operator pair instead of explicit, stored matrices. We modify the ASTRA Toolbox to incorporate operations involving the actual transpose matrices, $A^\top$ and $B^\top$, required by the AB- and BA-GKB algorithms.

\begin{figure}[ht!]
 \centering
  {\includegraphics[width=.6\textwidth]{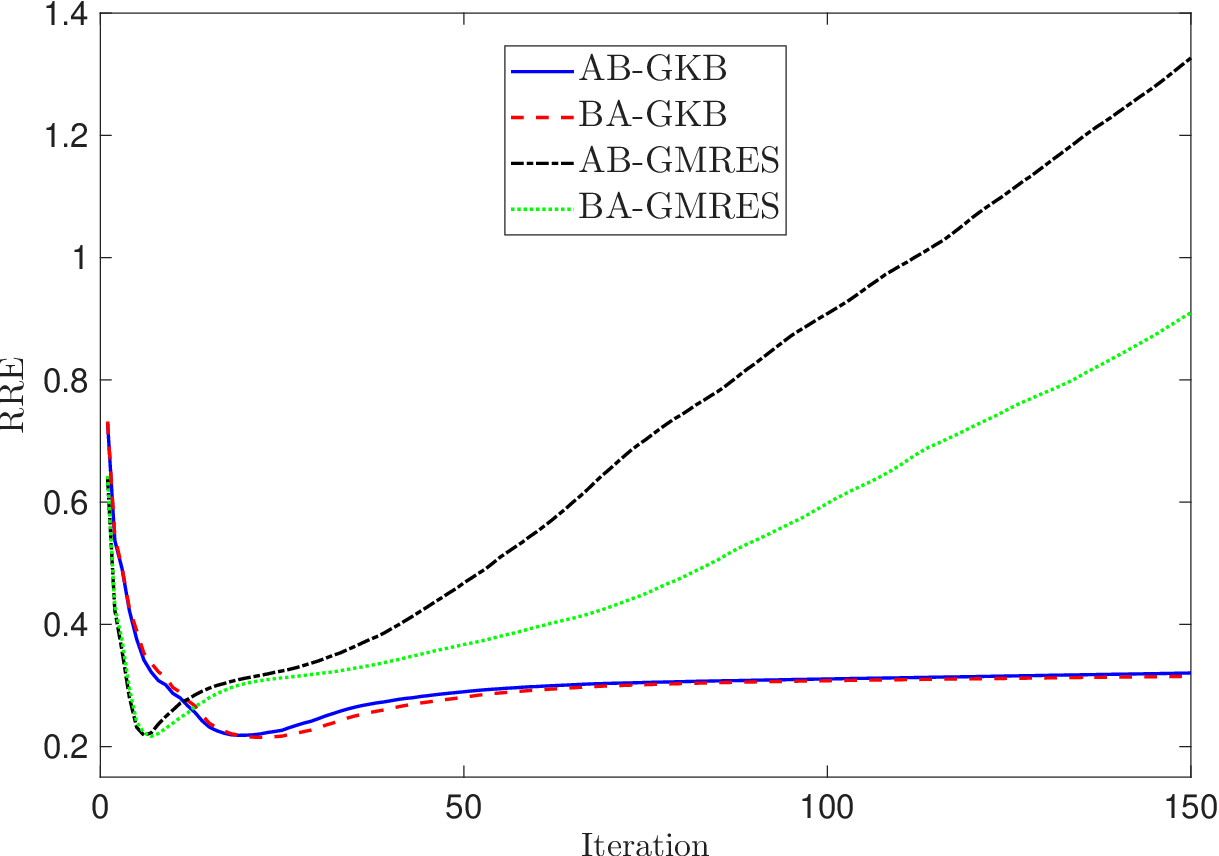}} 
 \caption{RRE for different methods for unmatched projectors  in \Cref{ex:problem1}.\label{Fig: mri RRE}}  
 \end{figure}

\begin{figure}[ht!]
  \centering
   \subfloat[AB-GKB]{\label{Fig:mri ABGKBDP}\includegraphics[width=.20\textwidth]{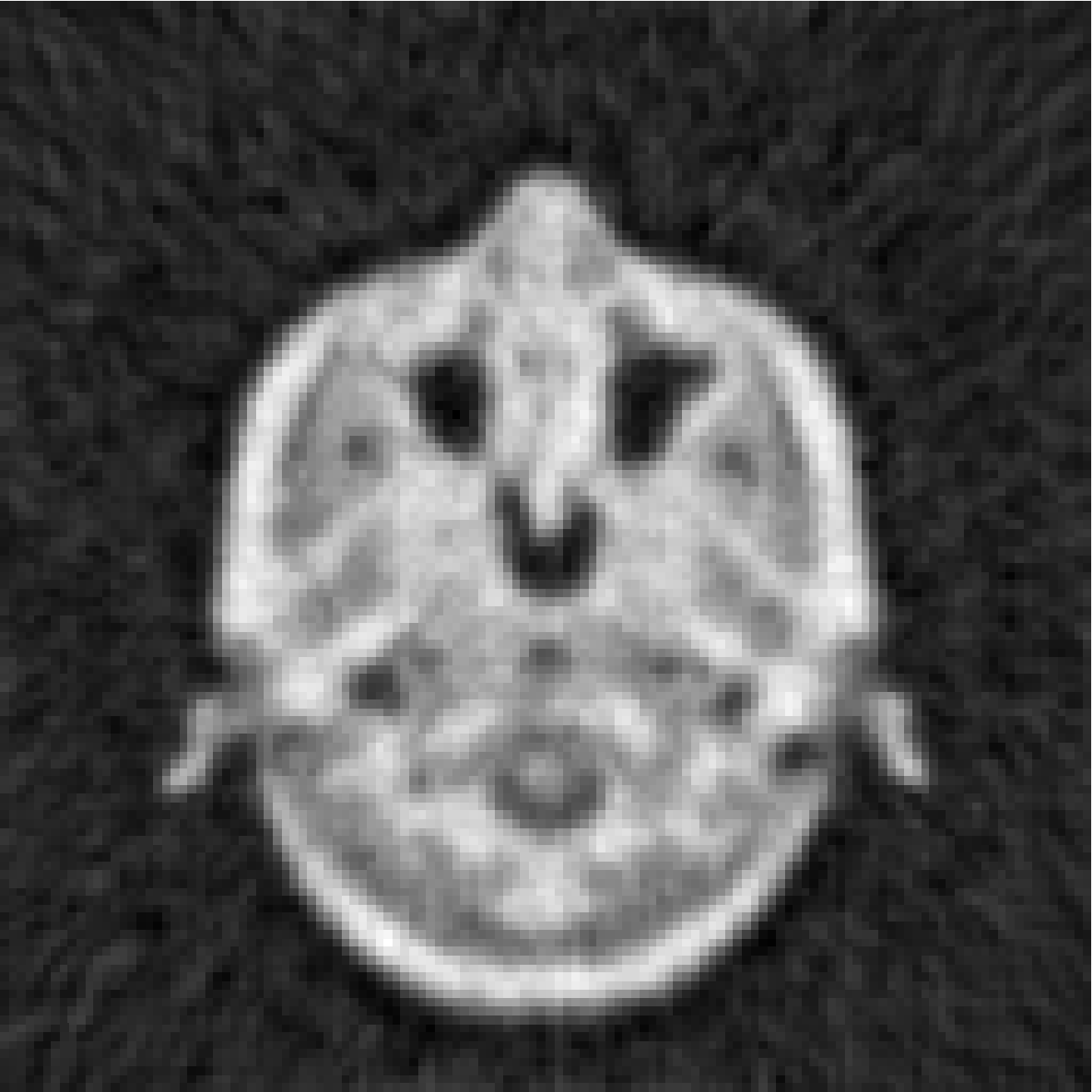}}   \ 
     \subfloat[BA-GKB]{\label{Fig:mri BAGKBDP}\includegraphics[width=.20\textwidth]{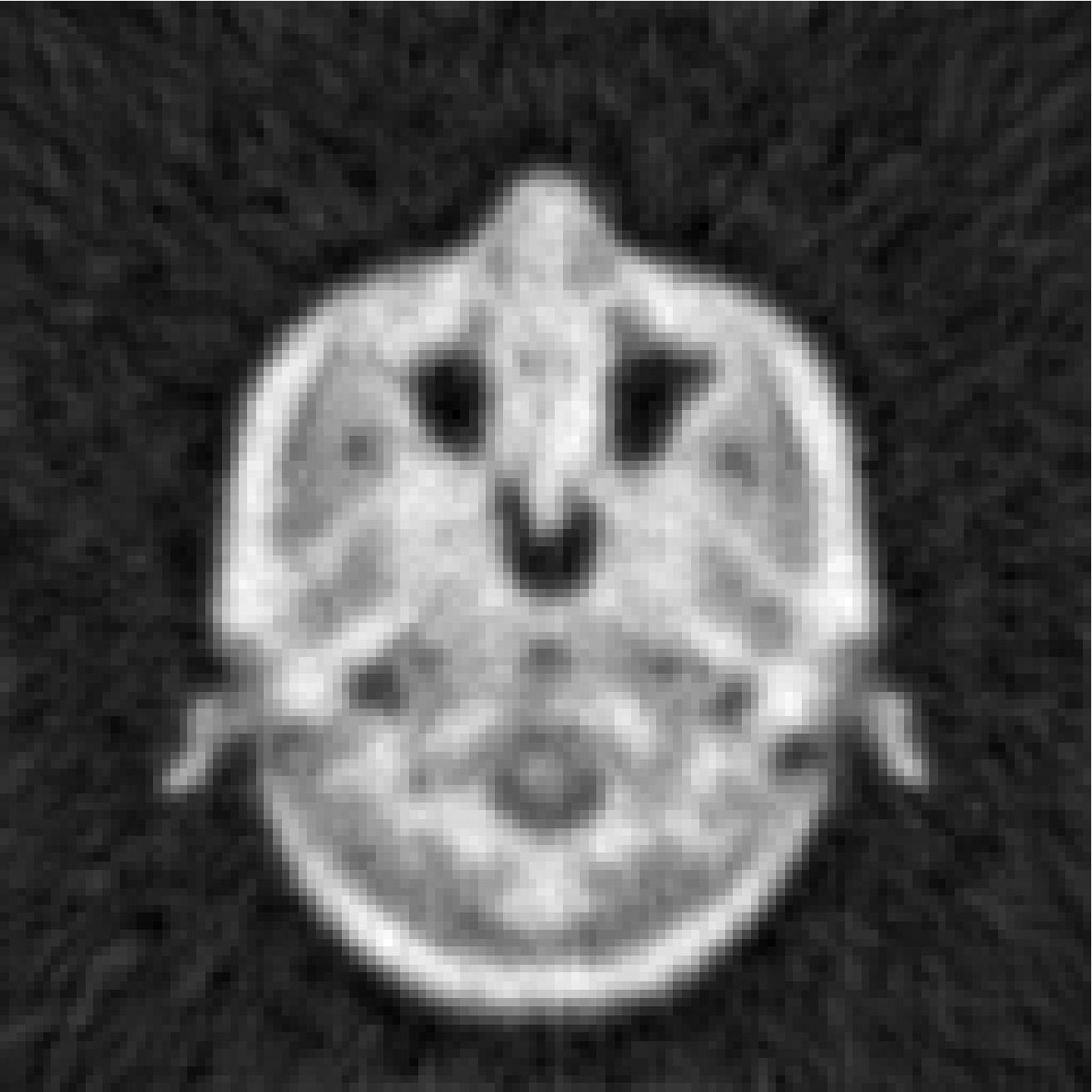}} \
  \subfloat[AB-GMRES]{\label{Fig:mri ABGNRESDP}\includegraphics[width=.20\textwidth]{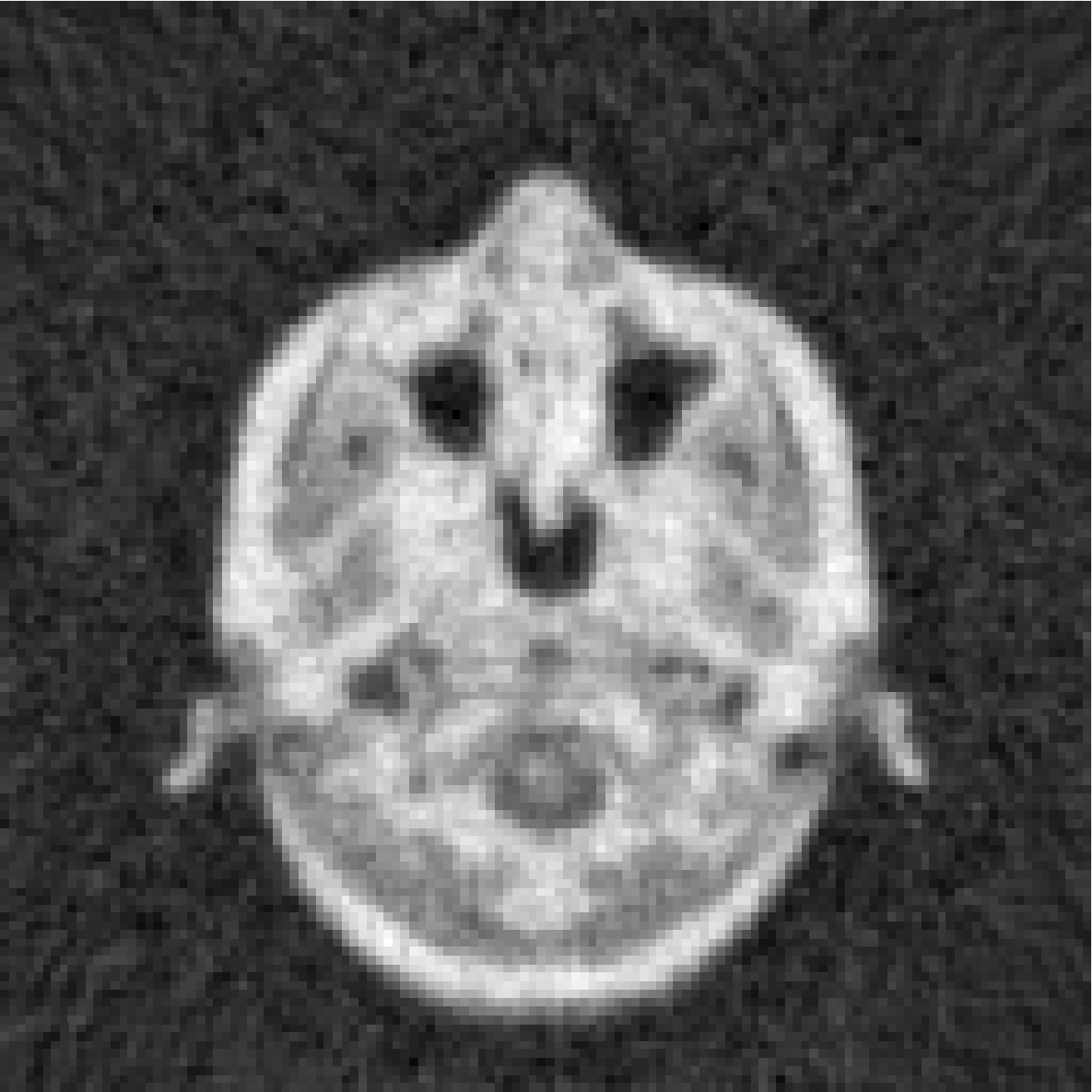}}  \
  \subfloat[BA-GMRES]{\label{Fig:mri BAGMRES}\includegraphics[width=.20\textwidth]{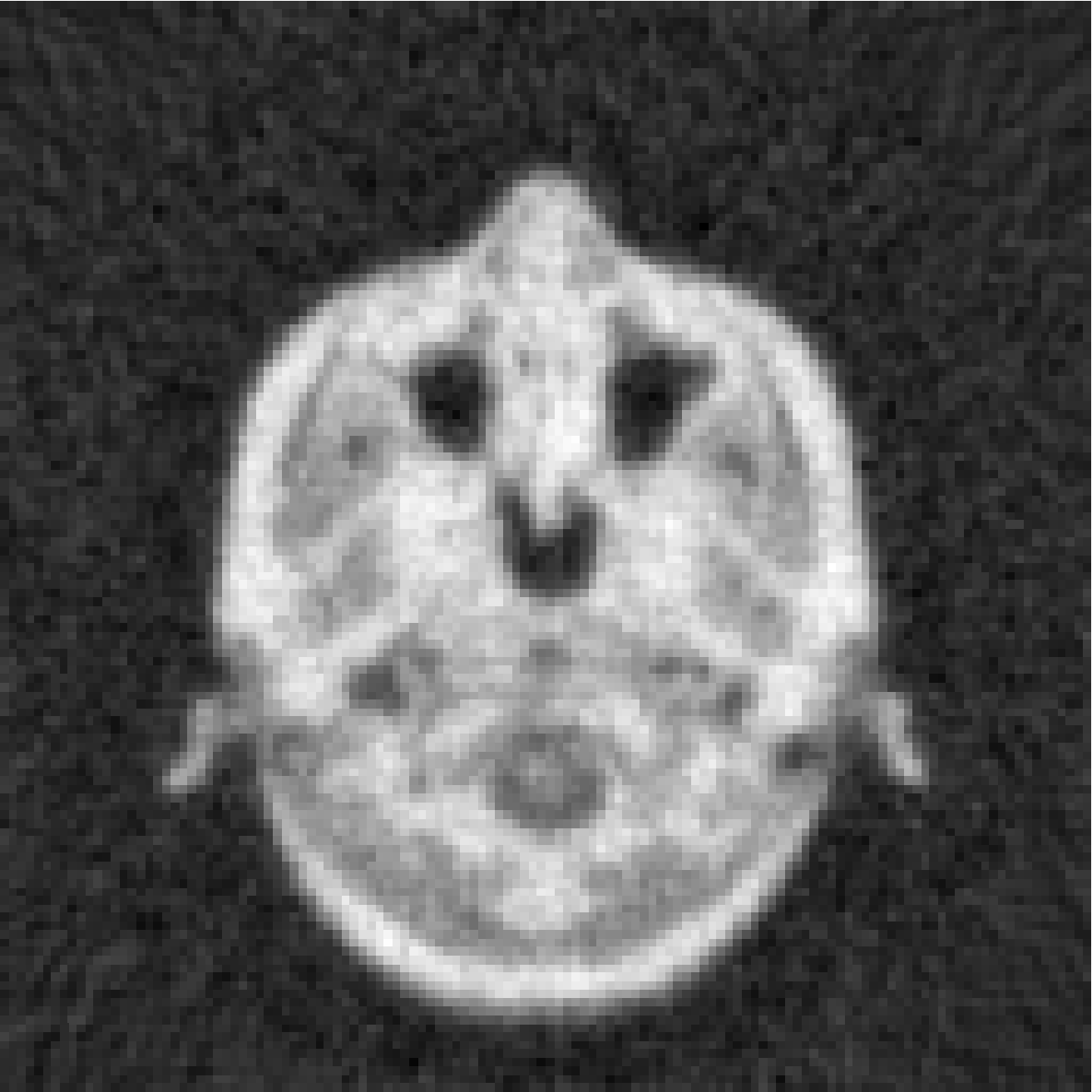}}  \\
    \subfloat[AB-GKB]{\label{Fig:mri ABGKBRNS}\includegraphics[width=.20\textwidth]{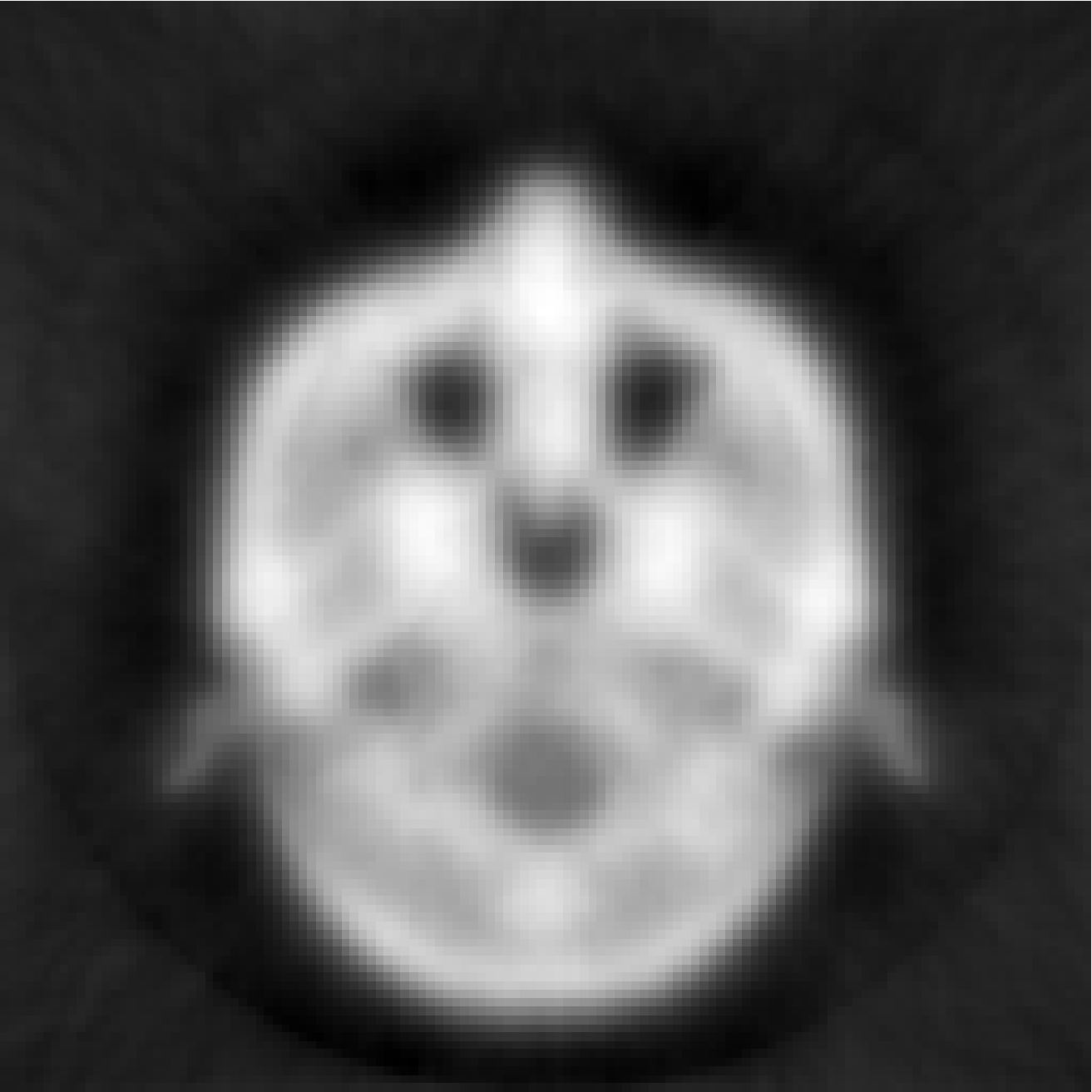}}   \ 
     \subfloat[BA-GKB]{\label{Fig:mri BAGKBRNS}\includegraphics[width=.20\textwidth]{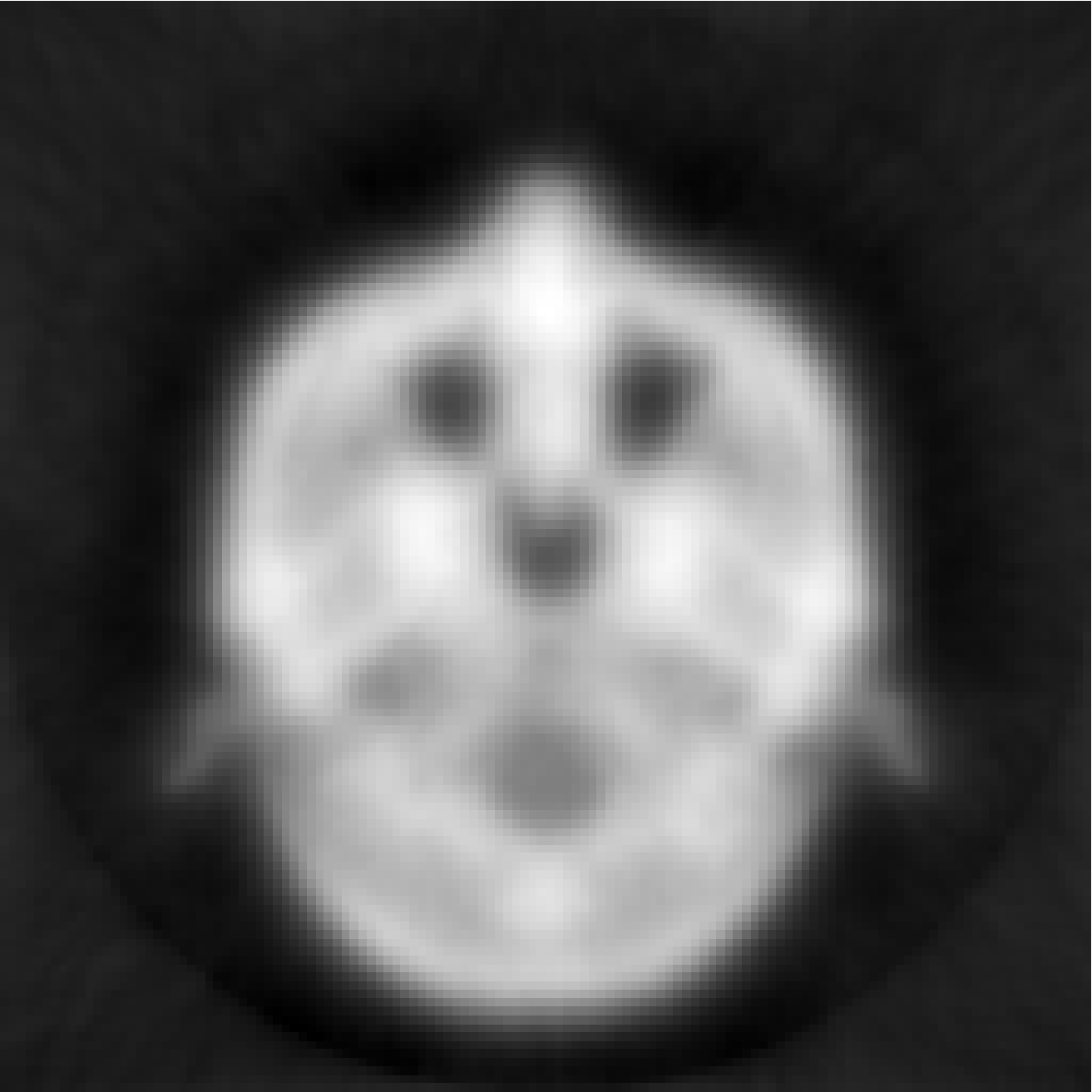}} \
  \subfloat[AB-GMRES]{\label{Fig:mri ABGNRESRNS}\includegraphics[width=.20\textwidth]{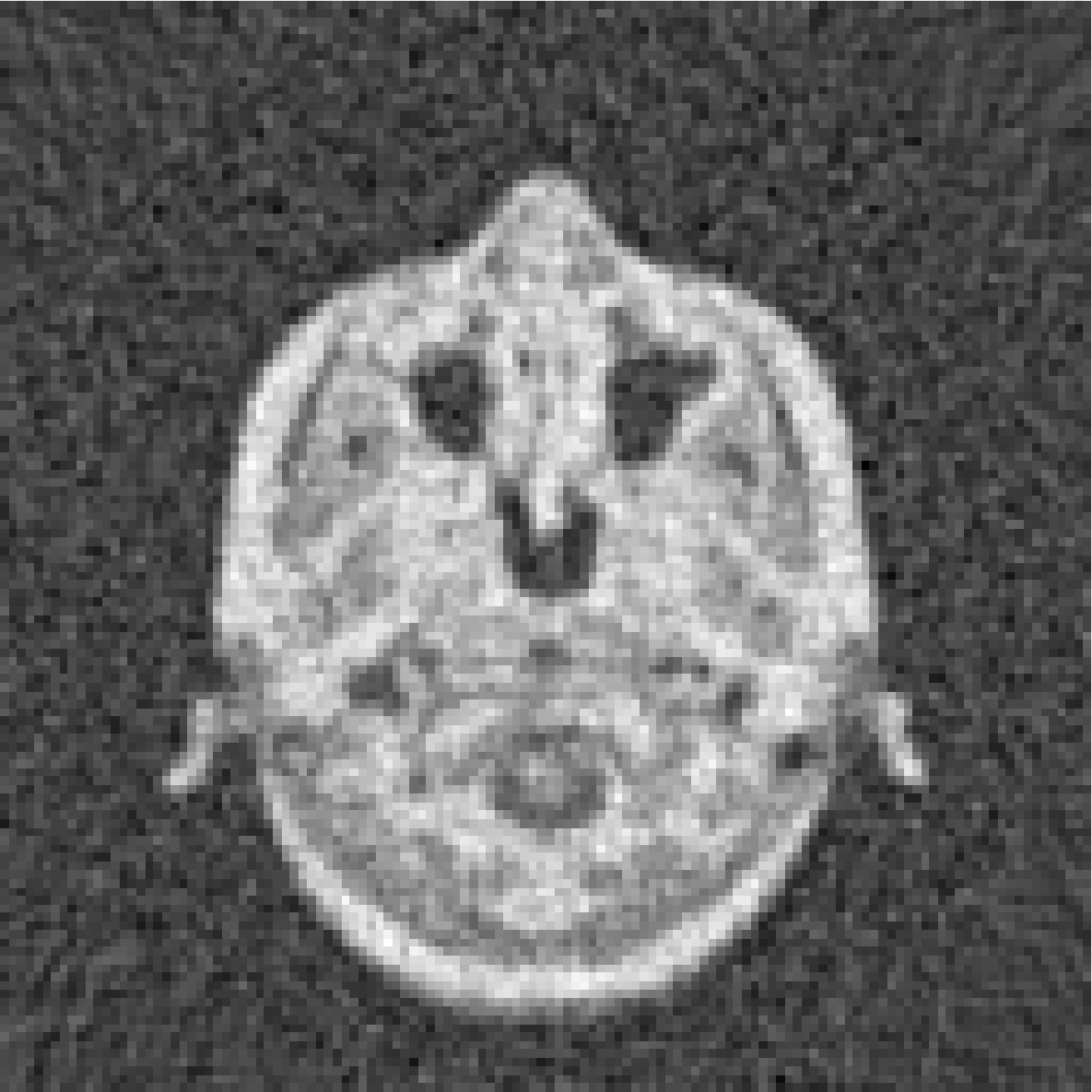}}  \
  \subfloat[BA-GMRES]{\label{Fig:mri BAGMRESRNS}\includegraphics[width=.20\textwidth]{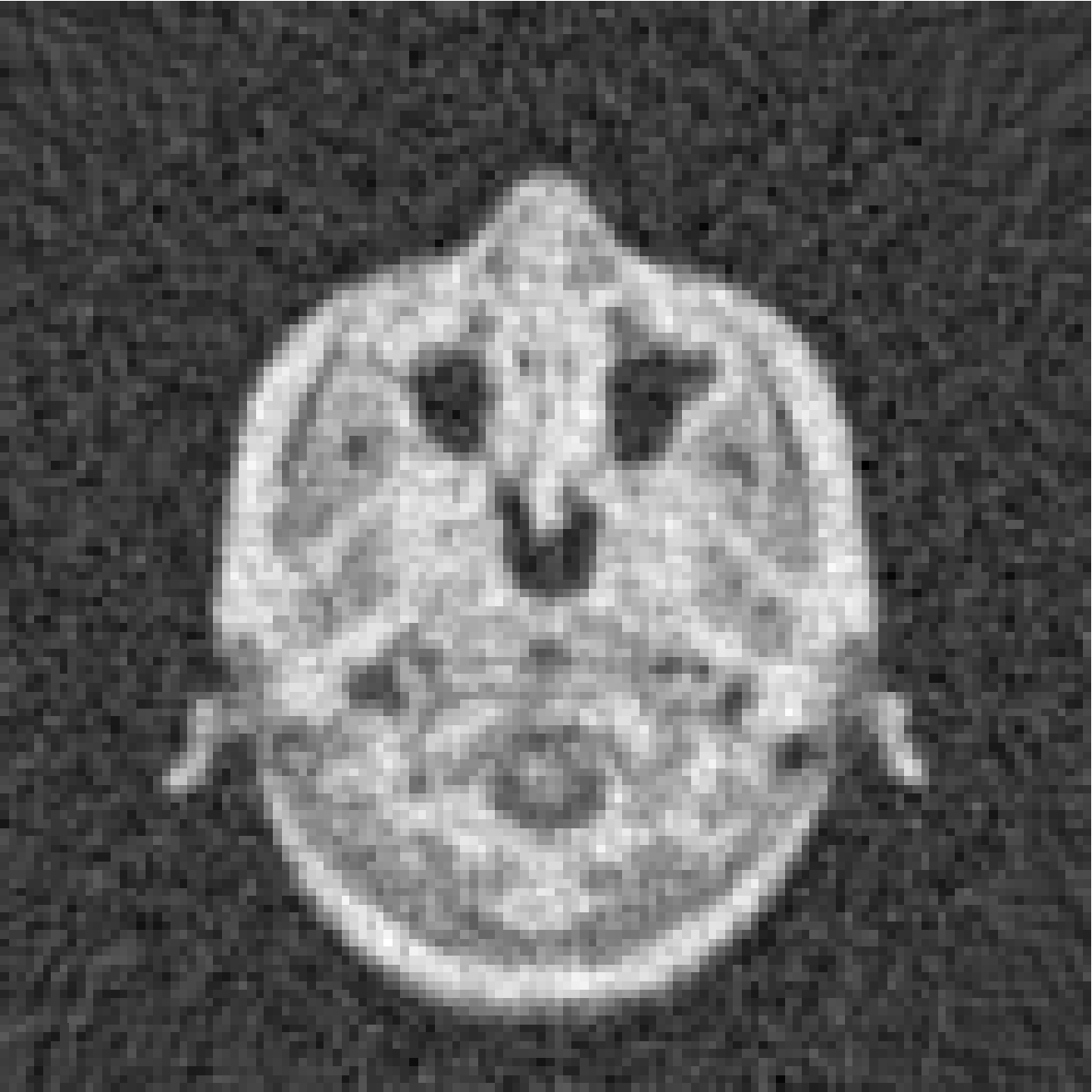}}  \

  \caption{Reconstructed images corresponding to \Cref{Fig:blurred mri}, obtained using DP and RNS stopping criteria in the first and second rows, respectively.}\label{Fig:Recons mri}  
  \end{figure}

From \Cref{Fig: mri RRE}, it is apparent that the AB- and BA-GKB methods exhibit less severe semiconvergence than their GMRES counterparts. Moreover, the GKB-based methods achieve lower RRE throughout the iterations, except during the first $11$ steps. While the GKB methods stabilize after approximately $40$ iterations, the GMRES methods diverge rapidly after about $7$ iterations. The corresponding reconstructions with different stopping criteria are provided in \Cref{Fig:Recons mri}.

%Despite these differences, the reconstructions obtained with DP, shown in \Cref{Fig:Recons mri}, display no significant visual discrepancies across the methods.

The results in column $3$ of \Cref{Tab:mri Problem} show that the DP yields lower RRE than the RNS criterion. Notably, DP consistently terminates all algorithms near their minimum RRE, demonstrating its effectiveness. These findings highlight the capability of GKB methods with DP to produce stable solutions without sacrificing GMRES-level accuracy.  The iteration counts in column $2$, consistent with \Cref{Fig: mri RRE}, indicate that GKB methods generally require more iterations to converge than their GMRES counterparts.

 \begin{table}[htb]
 \footnotesize 
 \caption{\label{Tab:mri Problem} The corresponding \text{RRE} \eqref{eq:RE} computed at $k_{DP}$ \eqref{eq:DP} and $k_{RNS}$ \eqref{eq:RNS}, by all methods, applied to the problem with $\text{SNR}\approx 28$ in \Cref{ex:problem1}. The lowest RRE results are shown in boldface. The results obtained using the RNS for the
corresponding methods are provided in parentheses.}
 \begin{tabular}[t]{lccccccccc } 
 \toprule 
Method&Iteration&$\text{RRE}(\bfx)$  \\ \midrule
 AB-GKB&$17(7)$&$\boldsymbol{0.21}(0.32)$\\
 BA-GKB&$17(7)$&$0.22(0.34)$\\  
 AB-GMRES&$6(9)$&$\boldsymbol{0.21}(0.25)$\\ 
 BA-GMRES&$6(9)$&$\boldsymbol{0.21}(0.23)$\\ 
\bottomrule
 \end{tabular}    
 \end{table} 

%%%%%%%%%%%%%%%%%%%%%%%%%%%%%%%%%%%%%%%%%%%%%%%%%%%%%%%%%%%%%%%%%%%%%%%%%%%%%%%%%%%

\begin{figure}[ht!]
  \centering
   \subfloat[True image]{\label{Fig:True Astra}\includegraphics[width=.25\textwidth]{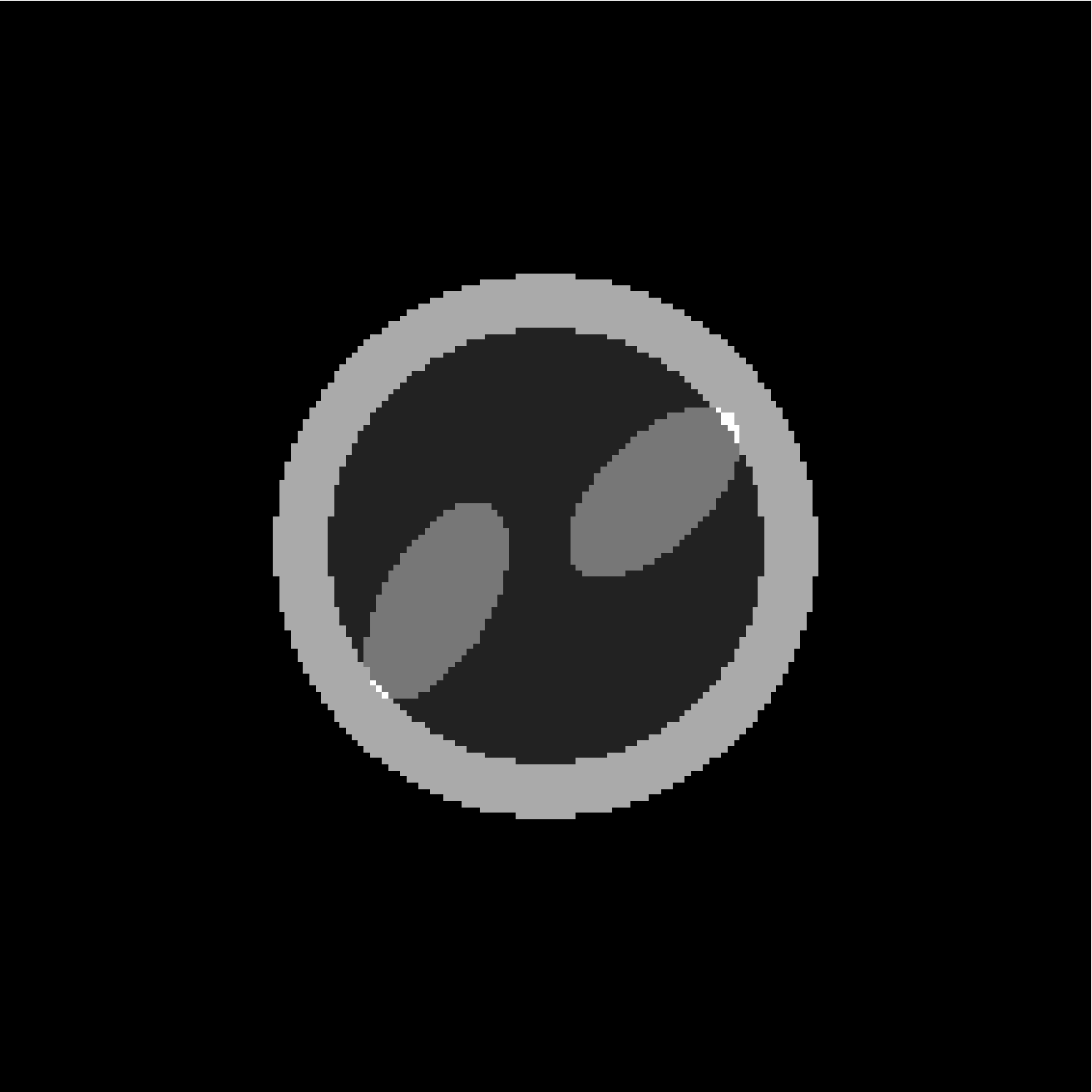}} \
   \subfloat[Noisy sinogram]{\label{Fig:blurred Astra}\includegraphics[width=.25\textwidth]{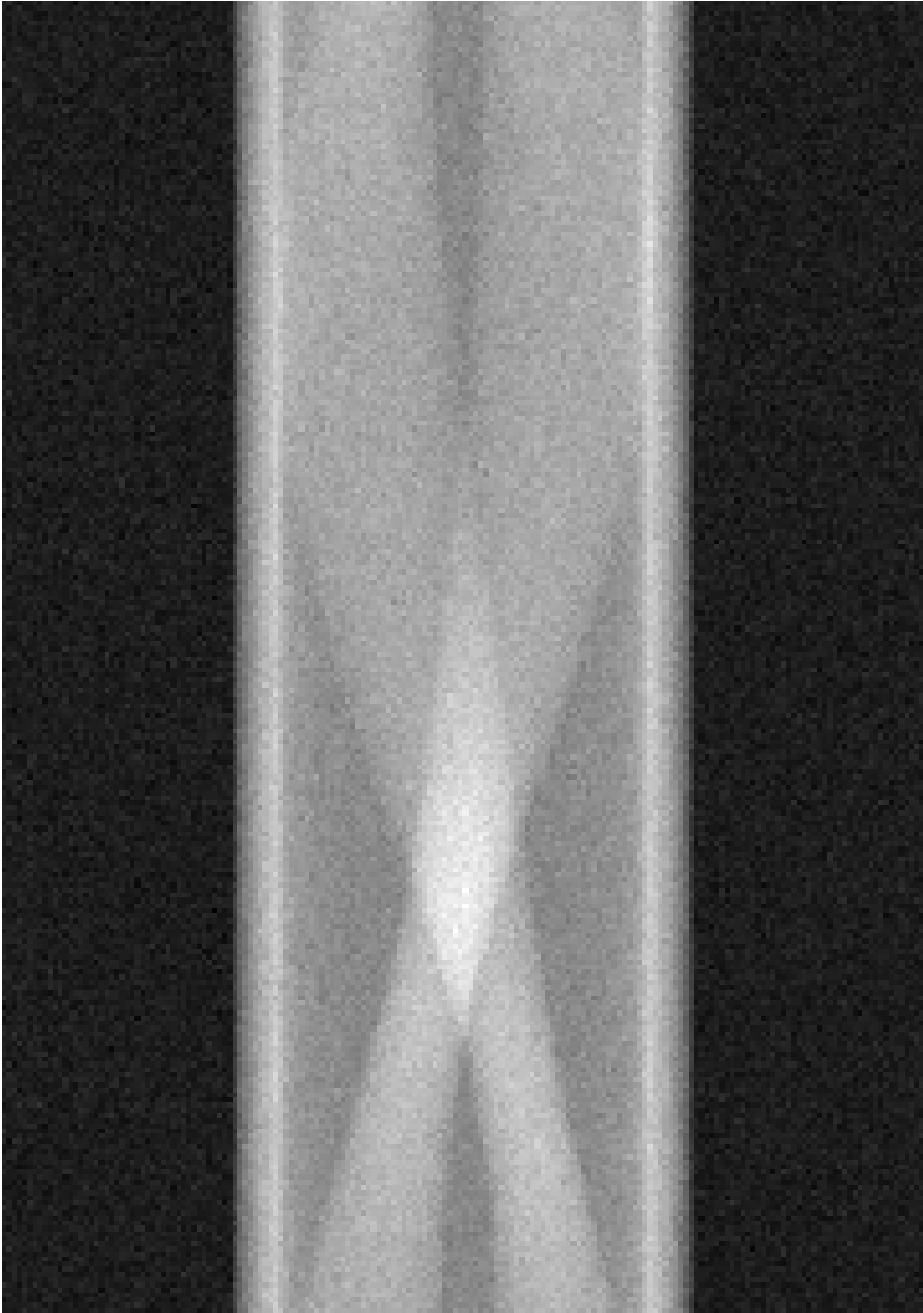}}
 \caption{The true image of $180 \times 180$ pixels and noisy sinogram with $7\%$ Gaussian noise, for the problem in \Cref{ex:problem2}. \label{Fig:Astra True and noisy}}  
  \end{figure}

\begin{example}\label{ex:problem2}
We consider a phantom image of size $180 \times 180$, shown in \Cref{Fig:True Astra}. The sinogram data are generated using parallel geometry through $256$ view angles. We use the ASTRA Toolbox \cite{van2015astra} to generate the forward projector $A \in \mathbb{R}^{46080\times 32400}$  and backprojector $B \in \mathbb{R}^{32400\times 46080}$.
Then we use $7 \%$ white Gaussian noise, corresponding to SNR of approximately $23$, to obtain the contaminated sinogram data $\bfb \in \mathbb{R}^{46080}$, shown in \Cref{Fig:blurred Astra}.
\end{example}

\begin{figure}[ht!]
 \centering
  {\includegraphics[width=.6\textwidth]{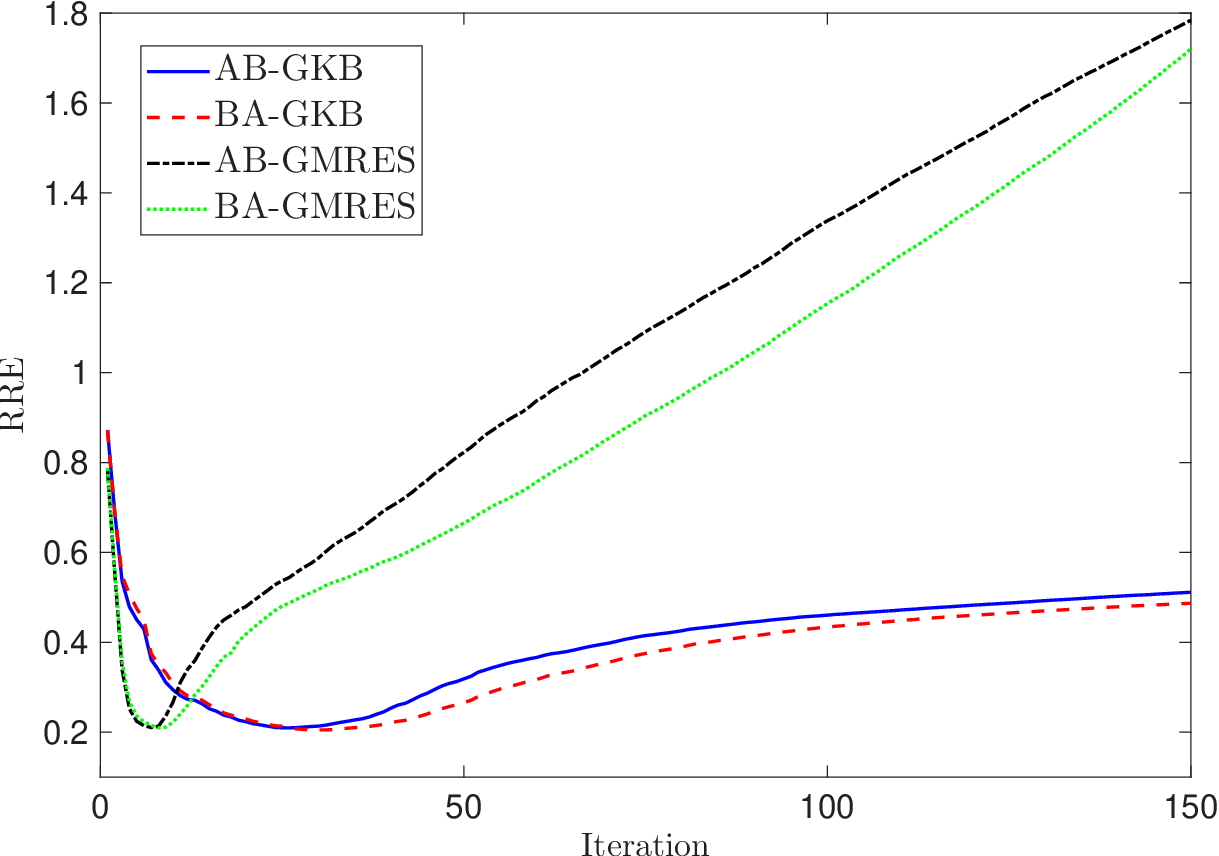}} 
 \caption{RRE for different methods for unmatched projectors  in \Cref{ex:problem2}.\label{Fig: Astra RRE}}  
 \end{figure}

\begin{figure}[ht!]
  \centering
   \subfloat[AB-GKB]{\label{Fig:ASTRA ABGKBDP}\includegraphics[width=.20\textwidth]{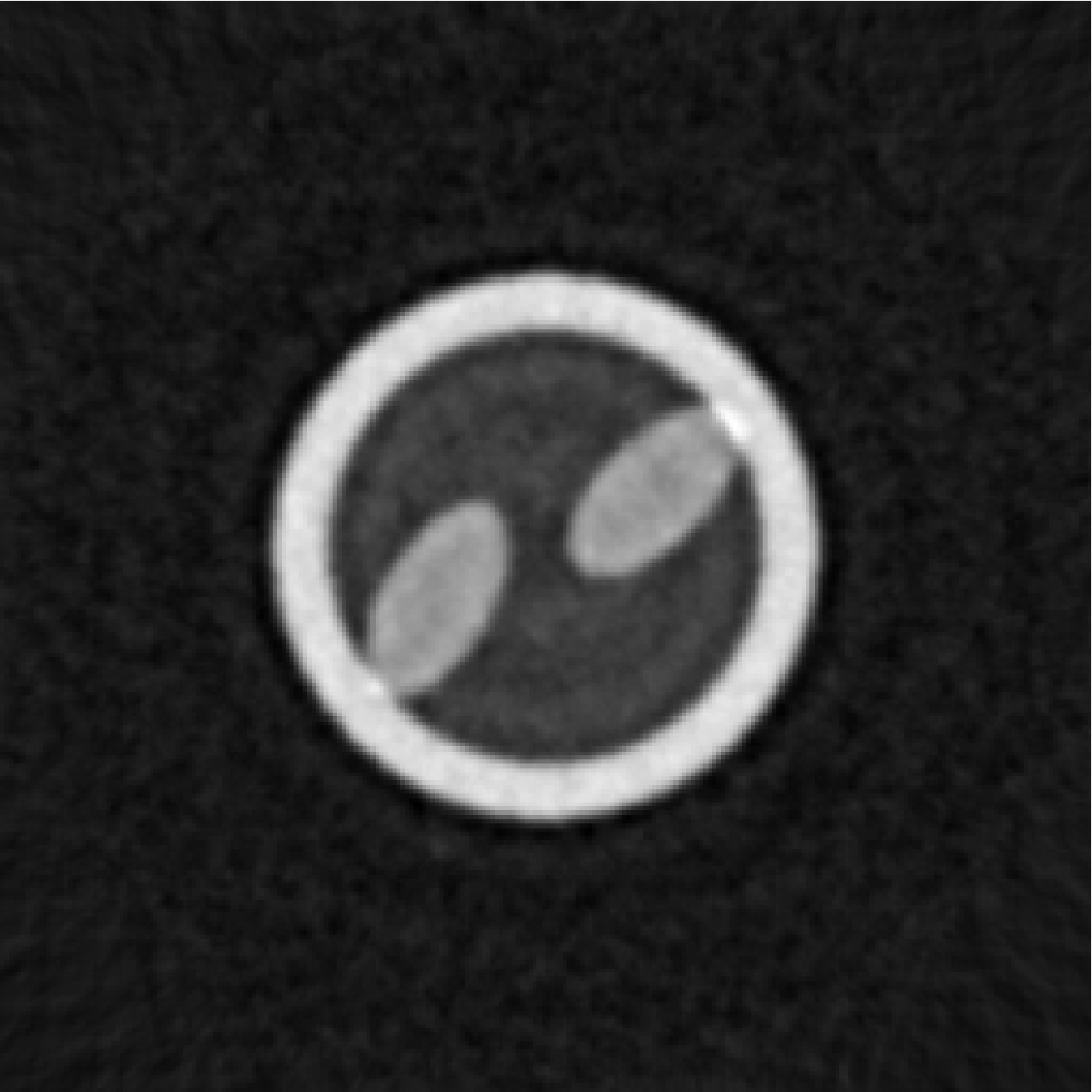}}   \ 
     \subfloat[BA-GKB]{\label{Fig:ASTRA BAGKBDP}\includegraphics[width=.20\textwidth]{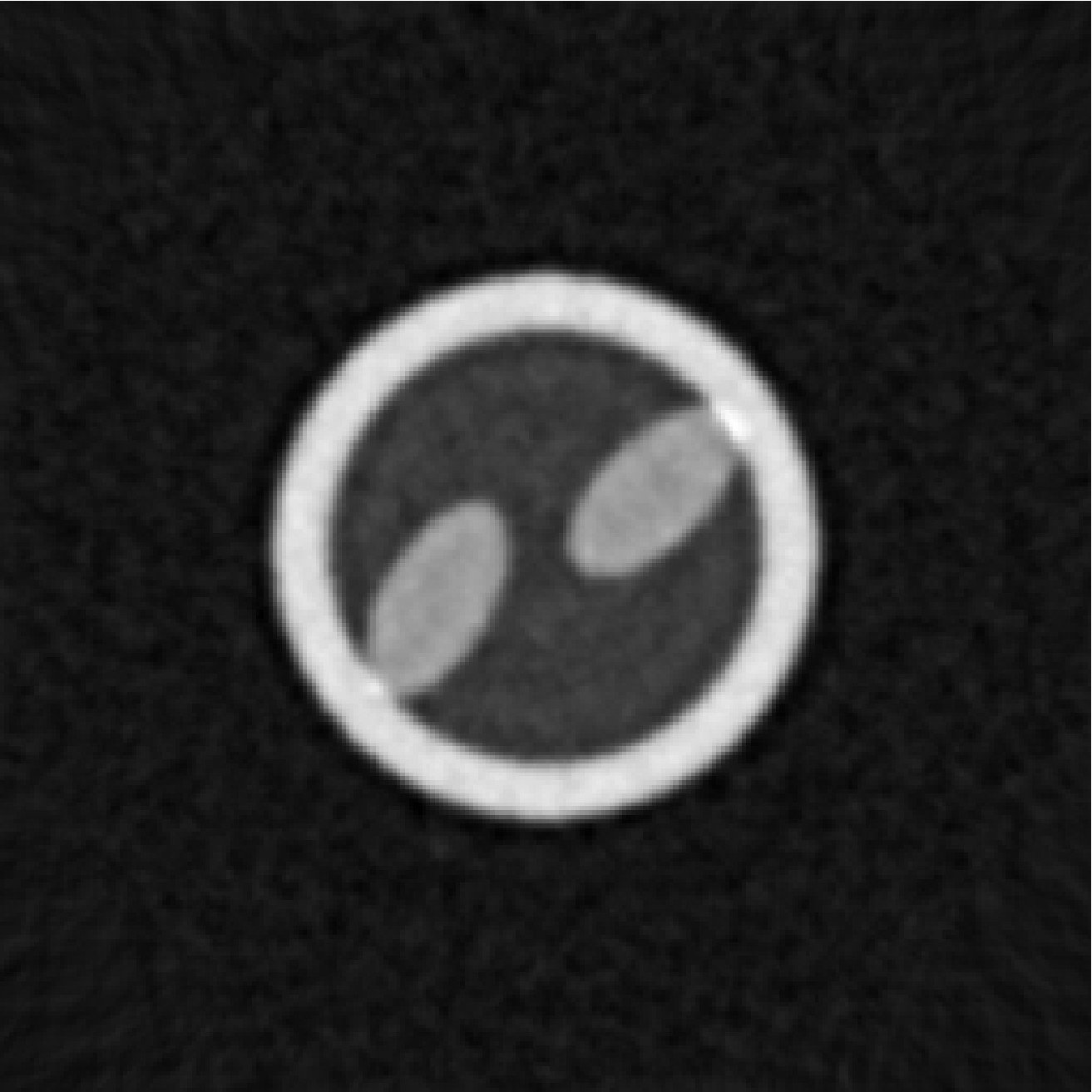}} \
  \subfloat[AB-GMRES]{\label{Fig:ASTRA ABGNRESDP}\includegraphics[width=.20\textwidth]{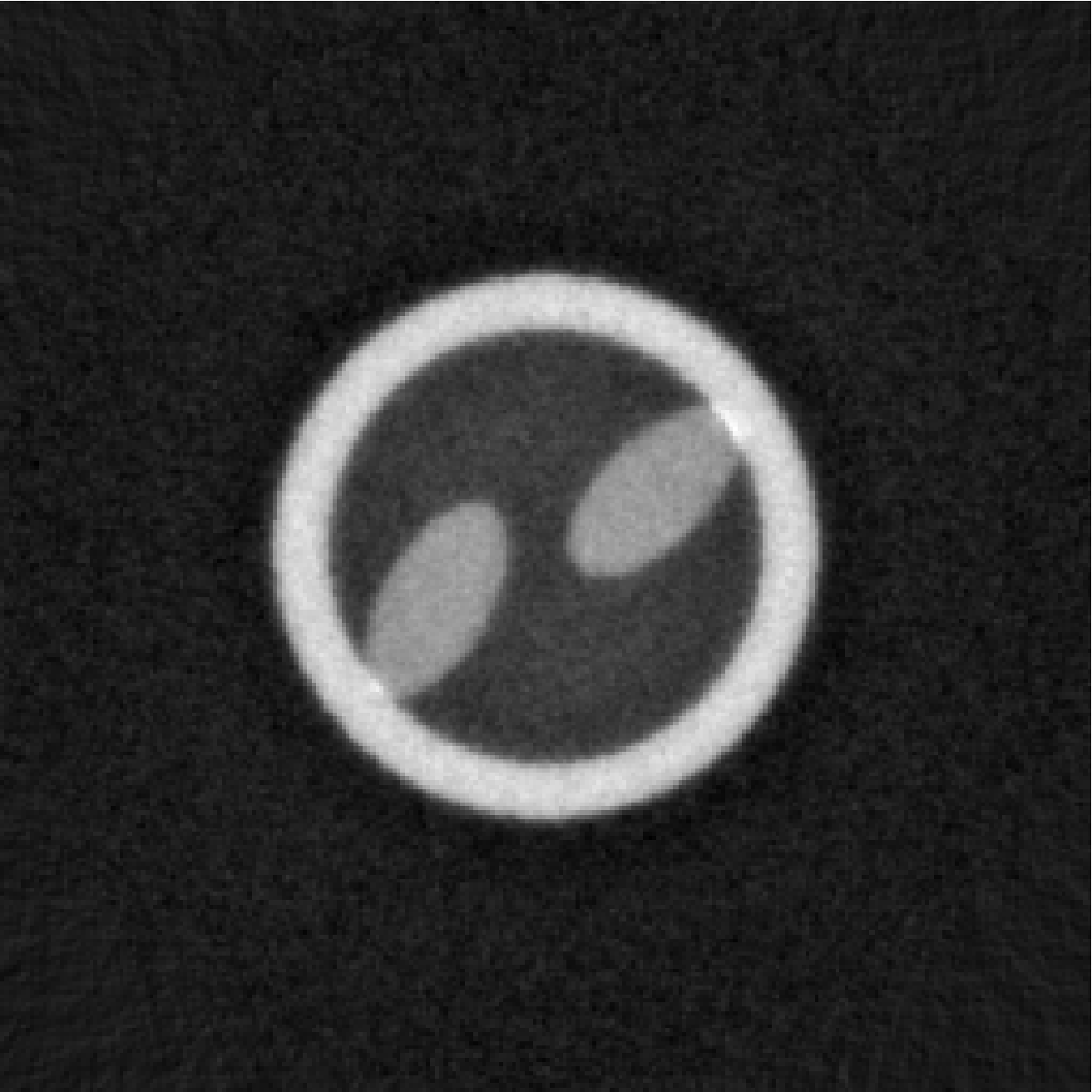}}  \
  \subfloat[BA-GMRES]{\label{Fig:ASTRA BAGMRES}\includegraphics[width=.20\textwidth]{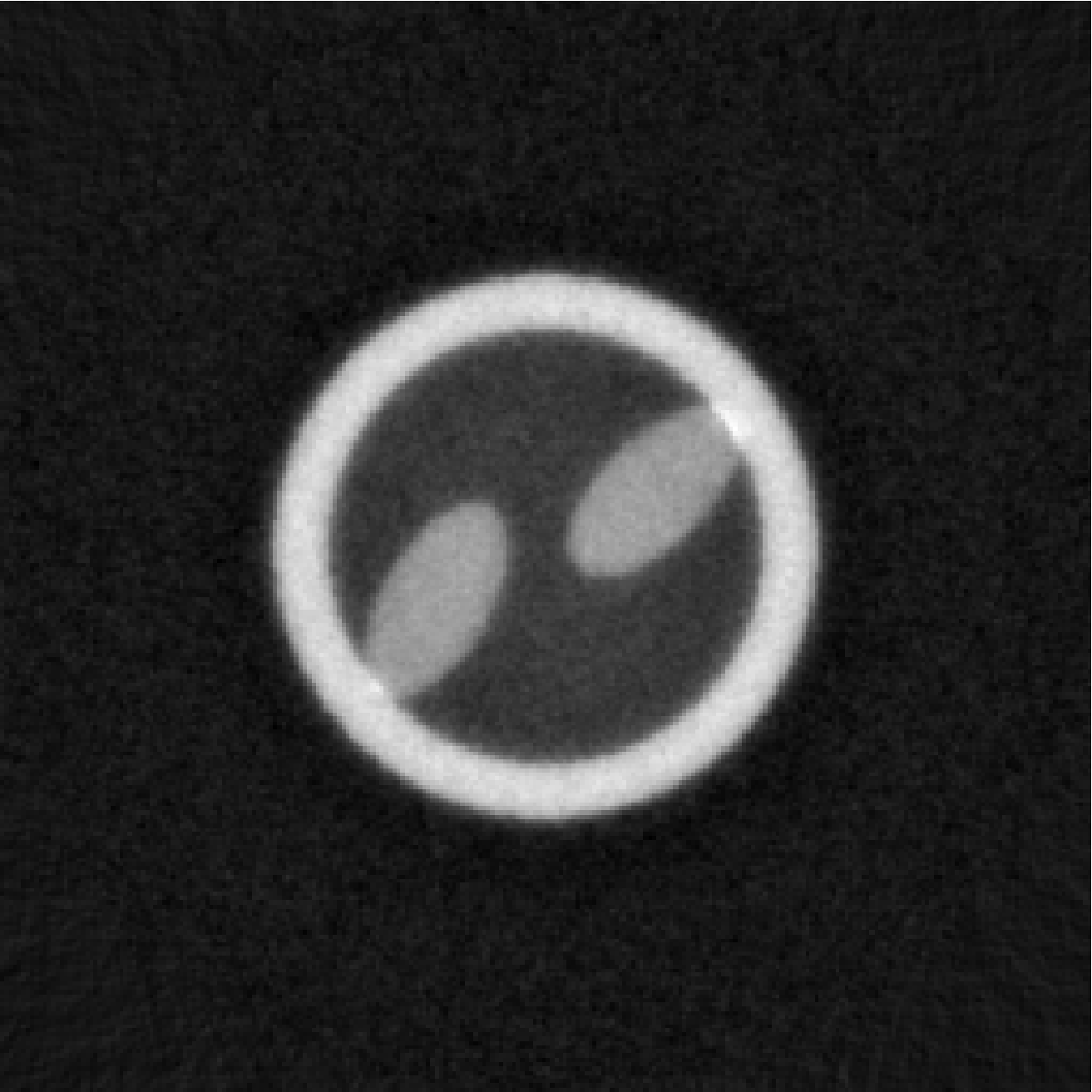}}  \\
   \subfloat[AB-GKB]{\label{Fig:ASTRA ABGKBRNS}\includegraphics[width=.20\textwidth]{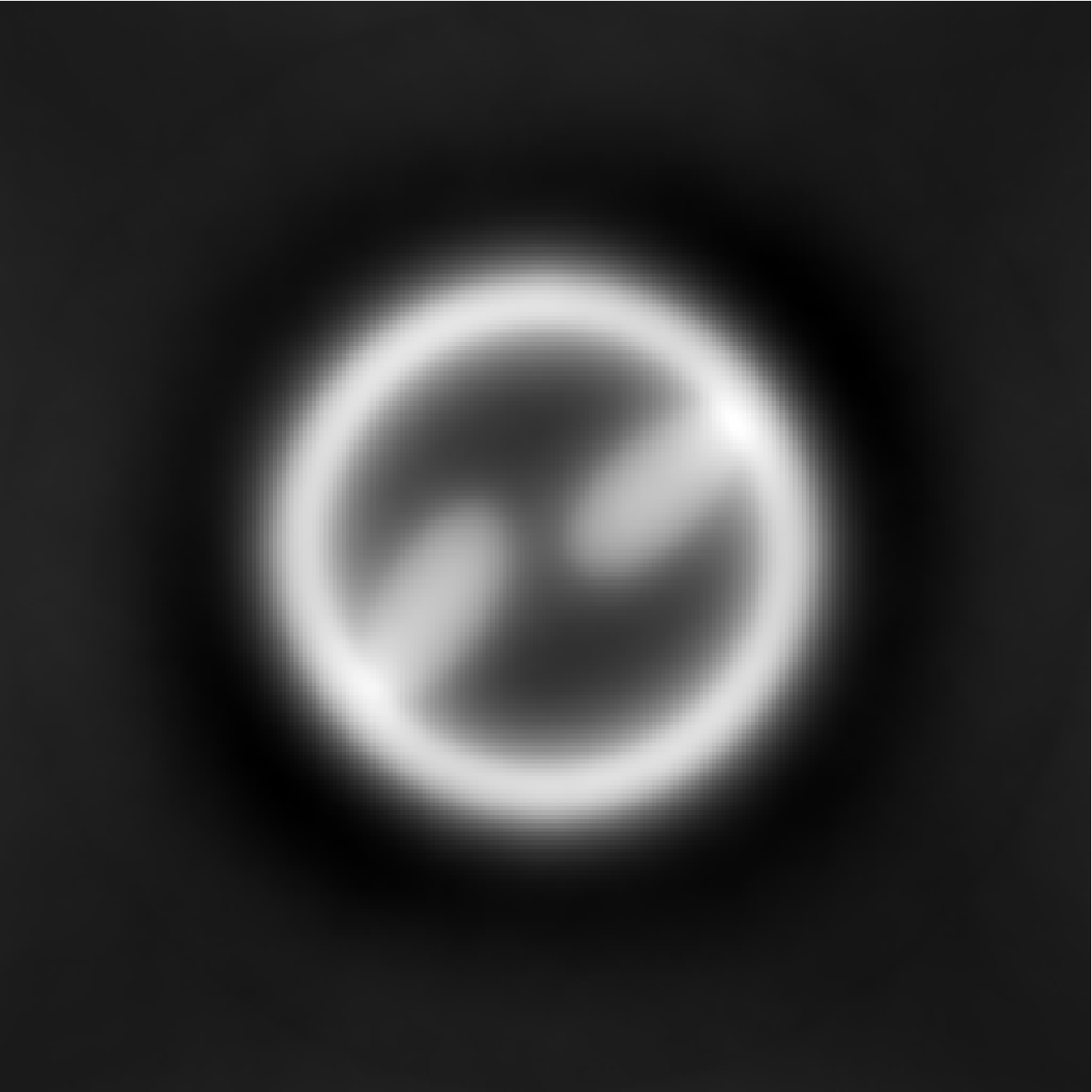}}   \ 
     \subfloat[BA-GKB]{\label{Fig:ASTRA BAGKBRNS}\includegraphics[width=.20\textwidth]{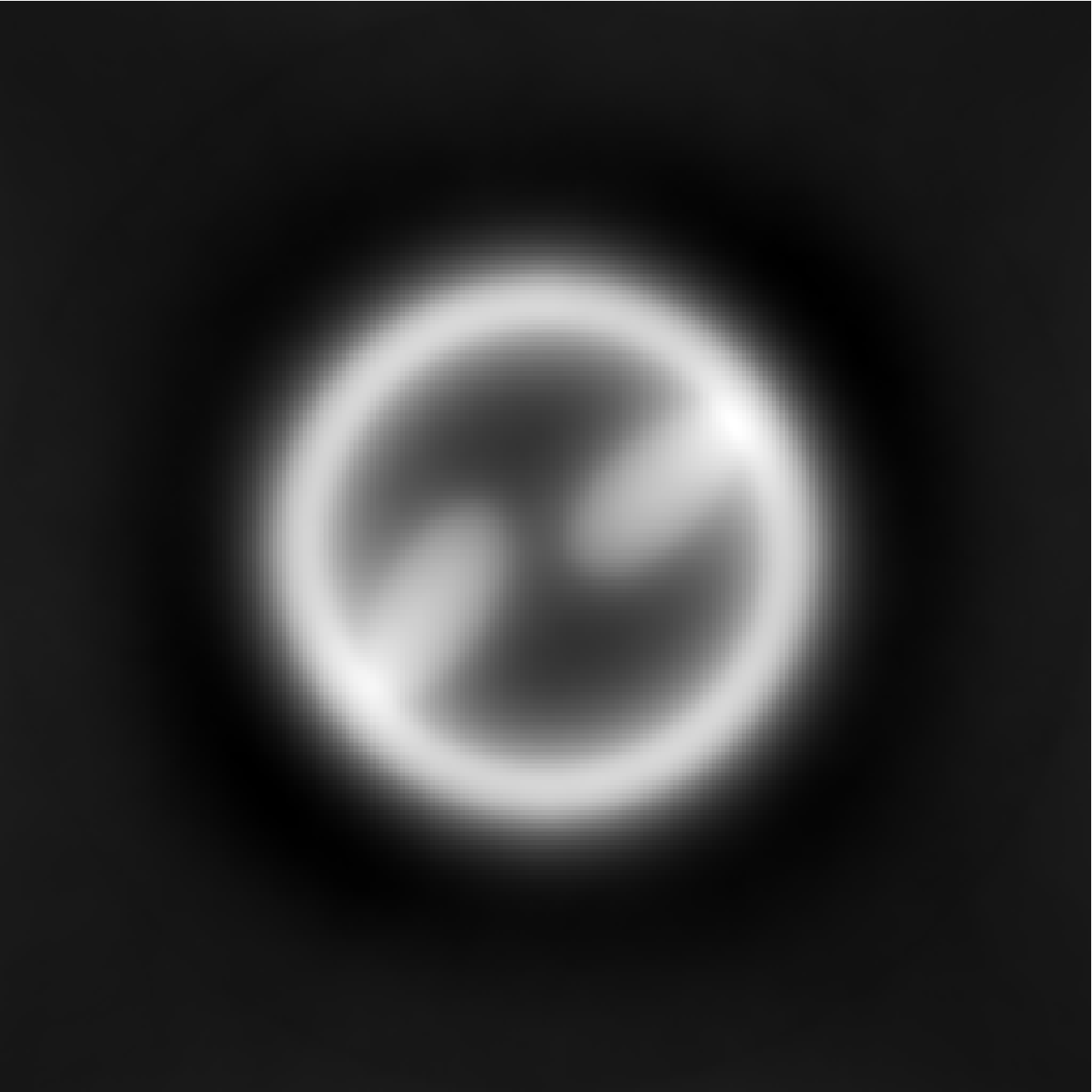}} \
  \subfloat[AB-GMRES]{\label{Fig:ASTRA ABGNRESRNS}\includegraphics[width=.20\textwidth]{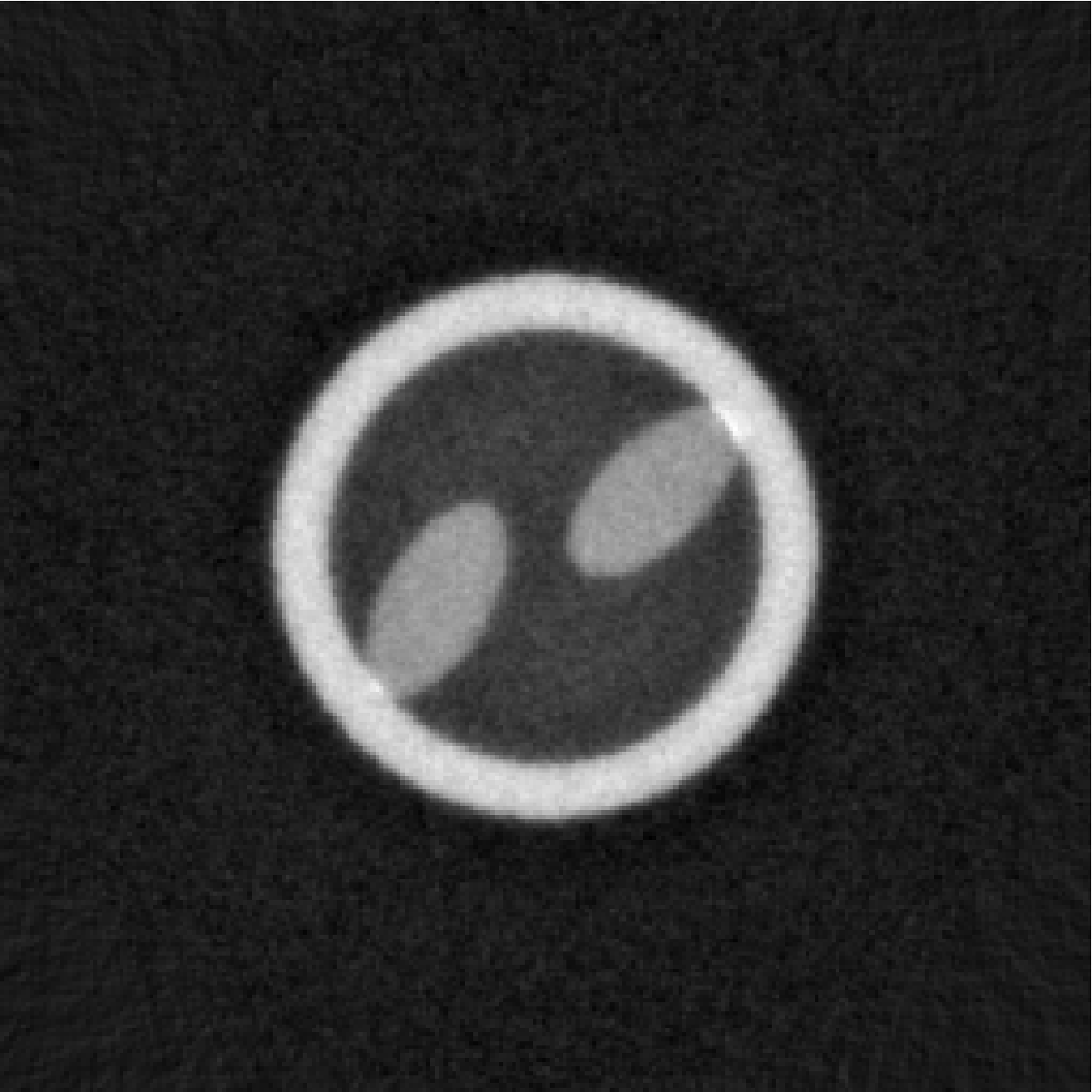}}  \
  \subfloat[BA-GMRES]{\label{Fig:ASTRA BAGMRESRNS}\includegraphics[width=.20\textwidth]{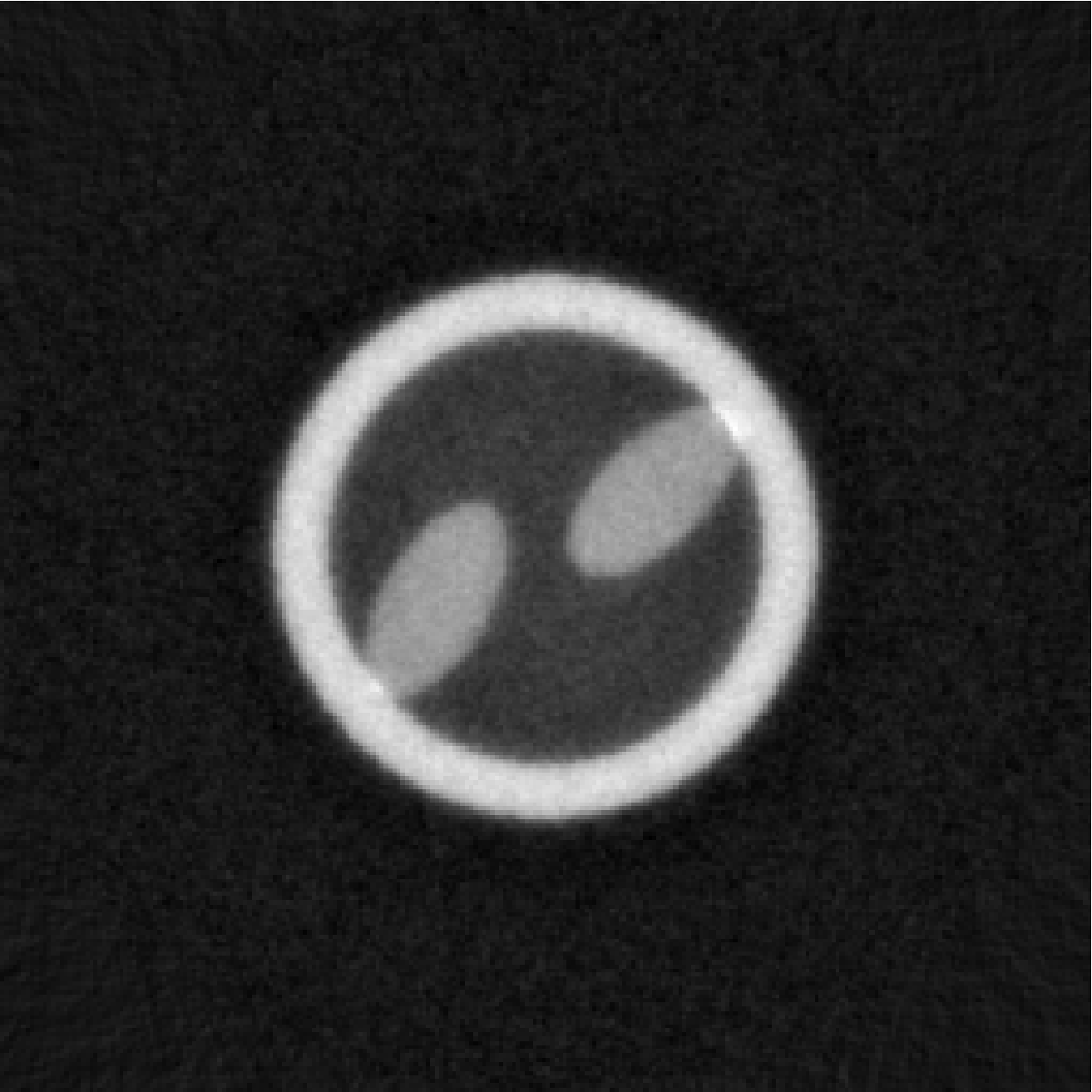}}  \

  \caption{Reconstructed images corresponding to \Cref{Fig:blurred Astra}, obtained using  DP and RNS stopping criteria in the first and second rows, respectively.}\label{Fig:Recons Astra}  
  \end{figure}

 This slightly larger problem is used to assess the performance of the proposed algorithms in the context of an overdetermined system matrix and to evaluate computational complexity. It should be emphasized that GPU-accelerated operations provided by the ASTRA Toolbox are used to carry out matrix–vector multiplications involving $A^\top$ and $B^\top$. From \Cref{Fig: Astra RRE}, it is immediately apparent that the AB- and BA-GKB methods are less severe semiconvergence compared to the AB- and BA-GMRES methods. In terms of RRE, we can see that the AB- and BA-GKB algorithms outperform the AB- and BA-GMRES algorithms along the iterations, except in the first $10$ iterations, during which the GMRES-based approaches demonstrate slightly better performance. The corresponding reconstructions are presented in \Cref{Fig:Recons Astra}.

 The results in column $3$ of \Cref{Tab:Astra Problem} indicate that the DP achieves lower RRE than the RNS criterion for the AB- and BA-GKB methods. In contrast, DP and RNS yield comparable performance for the AB- and BA-GMRES methods, as the tolerance $\epsilon$ in RNS is explicitly tuned to minimize the RRE. While both criteria enhance GMRES performance, pairing GKB methods with DP guarantees stable solutions without sacrificing accuracy.

We compare the theoretical computational complexity discussed in \Cref{sec:comp costs} with the actual costs presented in column $4$ of \Cref{Tab:Astra Problem}. The results in \Cref{Tab:Astra Problem} confirm the validity of the computational analysis in \Cref{sec:comp costs}, demonstrating that the GMRES-based methods are faster than the GKB-based methods. Moreover, for the overdetermined coefficient matrix, the results indicate that the AB-GMRES and AB-GKB methods are slower than their BA counterparts. In this investigation, we do not consider memory access and usage.

 \begin{table}[htb]
 \footnotesize 
 \caption{\label{Tab:Astra Problem} The corresponding \text{RRE} \eqref{eq:RE} computed at $k_{DP}$ \eqref{eq:DP} and $k_{RNS}$ \eqref{eq:RNS}, by all methods, applied to the problem with $\text{SNR}\approx 23$ in \Cref{ex:problem2}. The lowest RRE results are shown in boldface. Timings are recorded for $150$ iterations of all algorithms. The results obtained using the RNS for the
corresponding methods are provided in parentheses.}
 \begin{tabular}[t]{lccccccccc } 
 \toprule 
Method&Iteration&$\text{RRE}(\bfx)$&Time(s)  \\ \midrule
 AB-GKB&$23(5)$&$\boldsymbol{0.21}(0.45)$&$3.63$\\
 BA-GKB&$24(5)$&$0.22(0.48)$&$2.76$\\  
 AB-GMRES&$6(6)$&$\boldsymbol{0.21}(\boldsymbol{0.21})$&$2.82$\\ 
 BA-GMRES&$7(6)$&$\boldsymbol{0.21}(0.22)$&$2.05$\\ 
\bottomrule
 \end{tabular}    
 \end{table} 

%%%%%%%%%%%%%%%%%%%%%%%%%%%%%%%%%%%%%%%%%%%%%%%%%%%%%%%%%%%%%%%%%%%%%%%%%%%%%%%%%%%

\begin{figure}[ht!]
  \centering
   \subfloat[True image]{\label{Fig:True Shepplogan}\includegraphics[width=.25\textwidth]{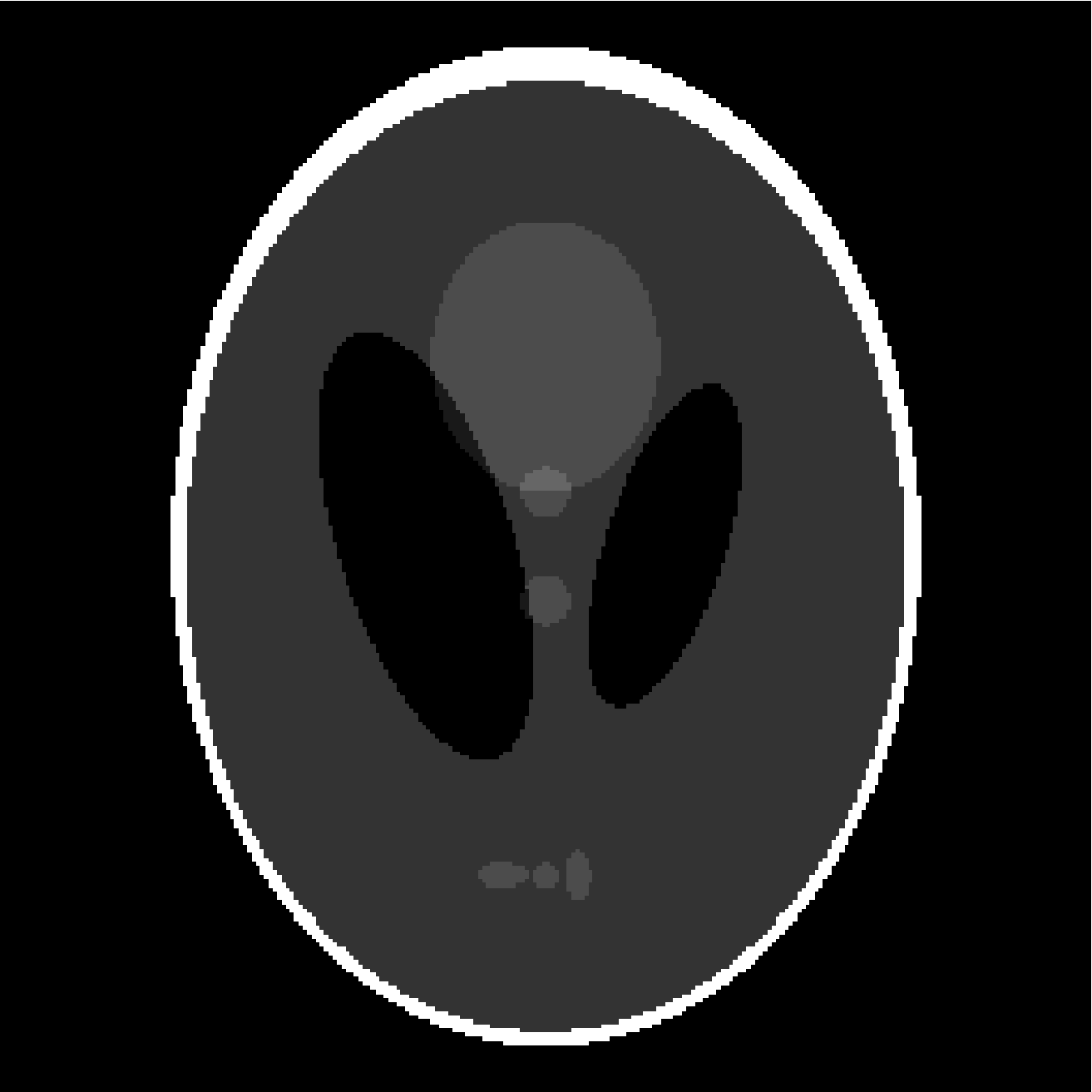}} \
   \subfloat[Noisy sinogram]{\label{Fig:blurred Shepplogan}\includegraphics[width=.25\textwidth]{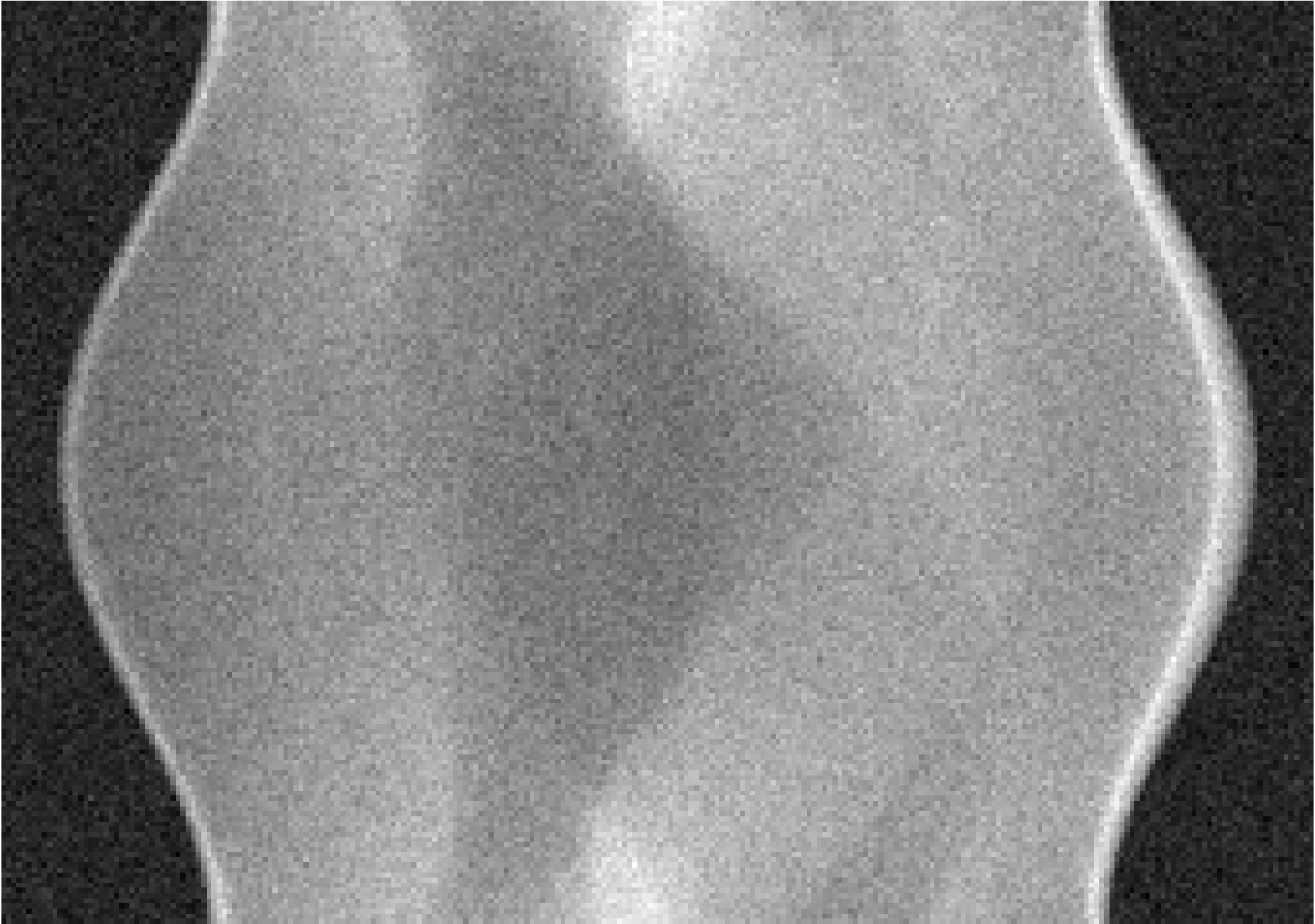}}
 \caption{The true image of $256 \times 256$ pixels and noisy sinogram with $10\%$ Gaussian noise, for the problem in \Cref{ex:problem3}. \label{Fig:Shepplogan True and noisy}}  
  \end{figure}

\begin{example}\label{ex:problem3}
 We consider the $\texttt{Shepplogan}$ image of size $256 \times 256$, shown in \Cref{Fig:True Shepplogan}. The sinogram data are generated using parallel geometry through $180$ view angles. We use the ASTRA Toolbox \cite{van2015astra} to generate the forward projector $A \in \mathbb{R}^{46080\times 65536}$  and backprojector $B \in \mathbb{R}^{65536\times 46080}$. The transpose matrices $A^\top$ and $B^\top$ are generated using a similar process. Then we use $10 \%$ white Gaussian noise, corresponding to SNR of approximately $20$, to obtain the contaminated sinogram data $\bfb \in \mathbb{R}^{65536}$, shown in \Cref{Fig:blurred Shepplogan}.
 \end{example}

\begin{figure}[ht!]
 \centering
  {\includegraphics[width=.6\textwidth]{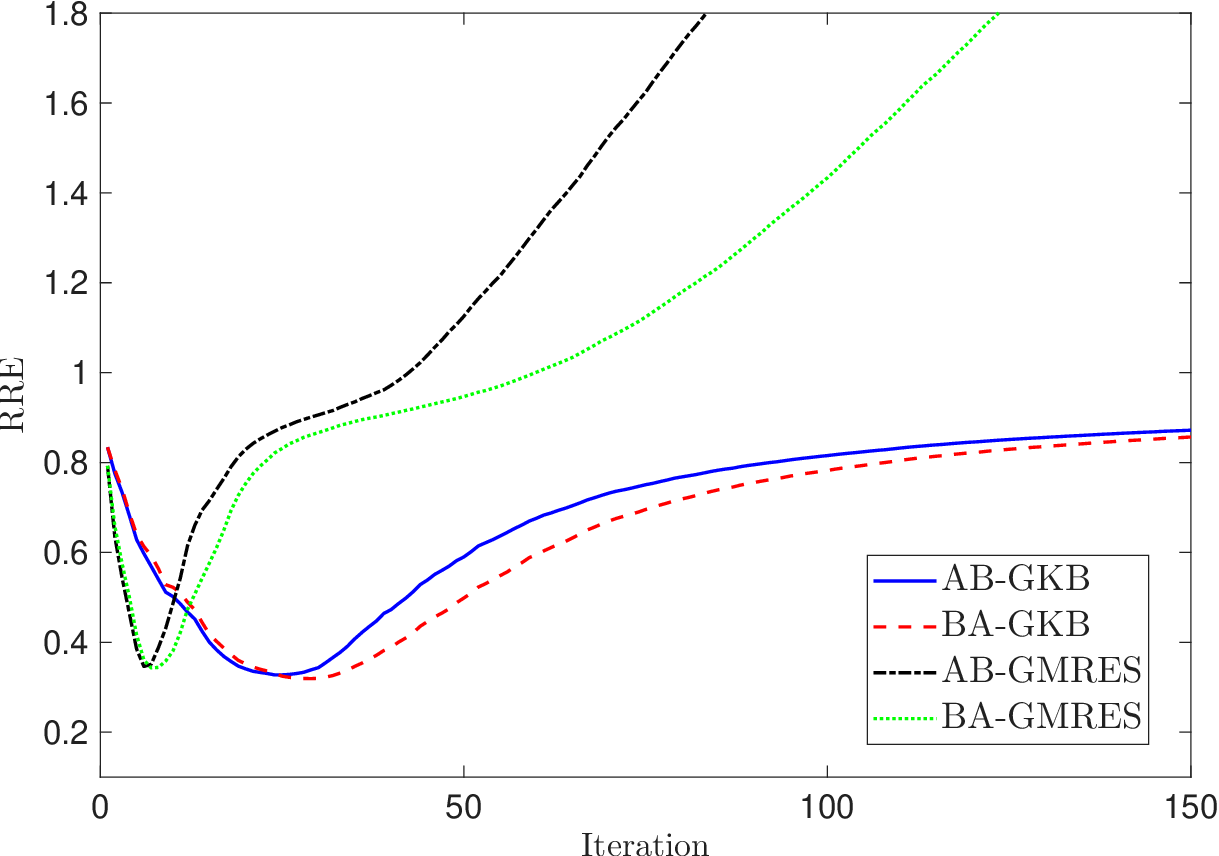}} 
 \caption{RRE for different methods for unmatched projectors in \Cref{ex:problem3}.\label{Fig:Shepplogan RRE}}  
 \end{figure}

\begin{figure}[ht!]
  \centering
   \subfloat[AB-GKB]{\label{Fig:Shepplogan ABGKBDP}\includegraphics[width=.20\textwidth]{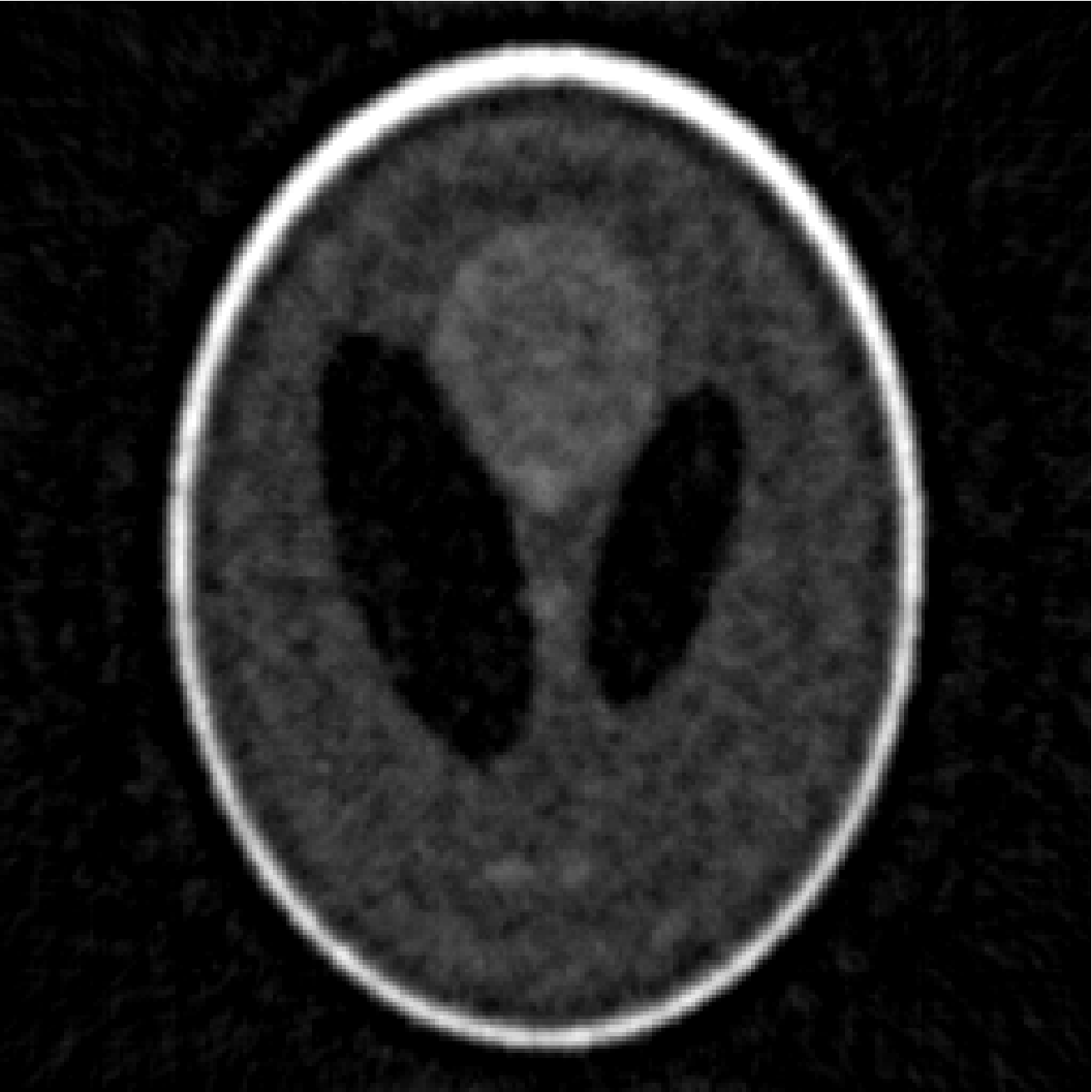}}   \ 
     \subfloat[BA-GKB]{\label{Fig:ShepploganShepplogan BAGKBDP}\includegraphics[width=.20\textwidth]{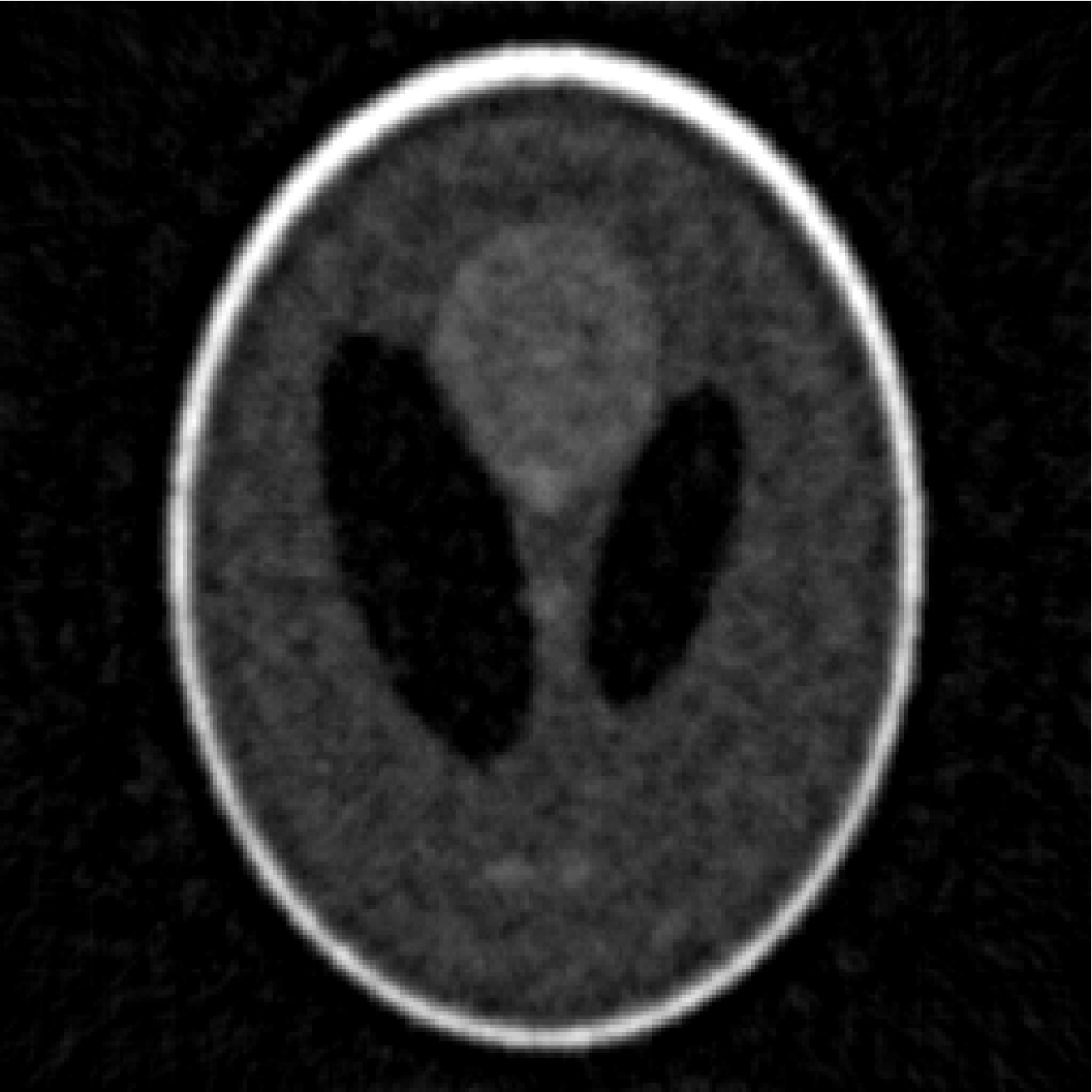}} \
  \subfloat[AB-GMRES]{\label{Fig:Shepplogan ABGNRESDP}\includegraphics[width=.20\textwidth]{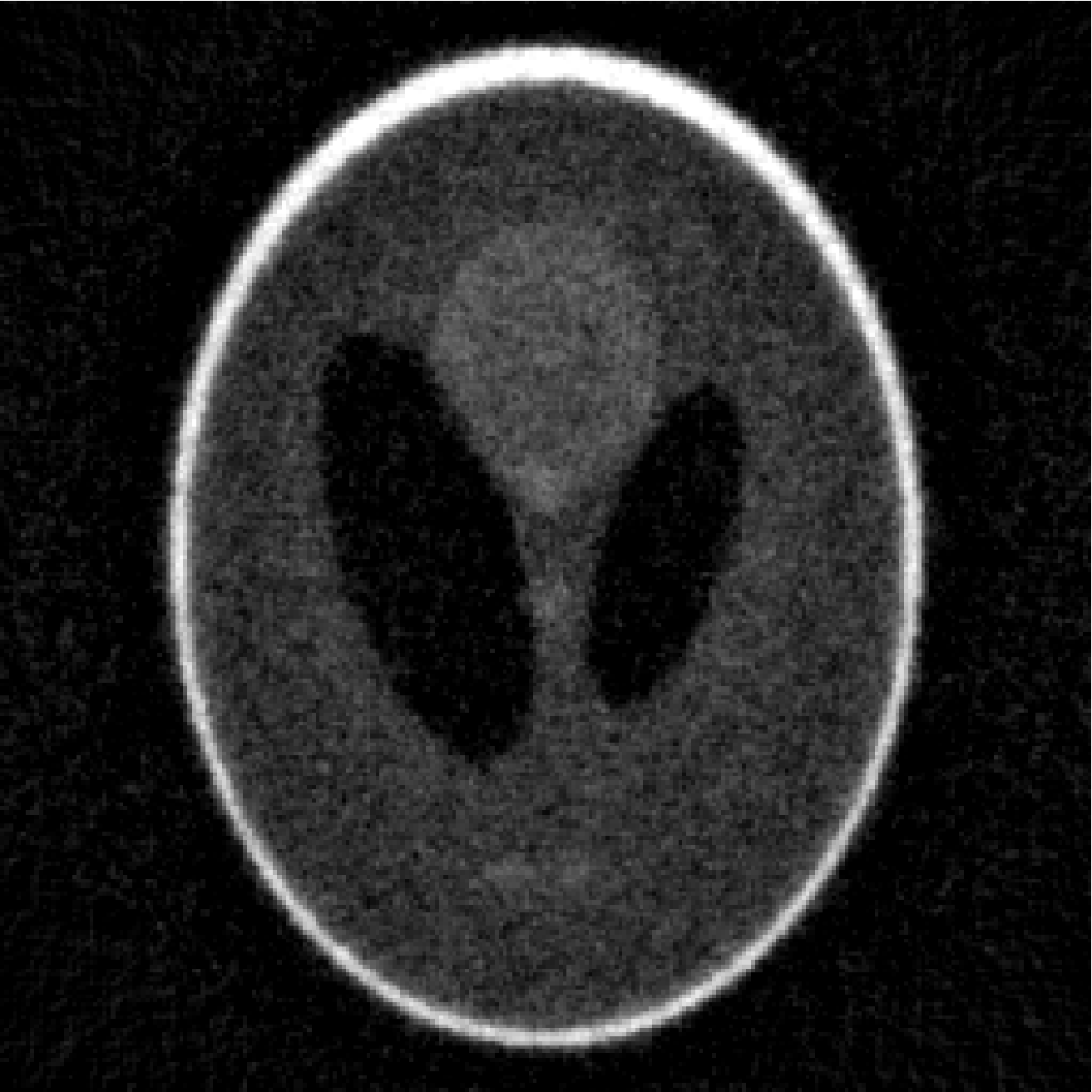}}  \
  \subfloat[BA-GMRES]{\label{Fig:Shepplogan BAGMRES}\includegraphics[width=.20\textwidth]{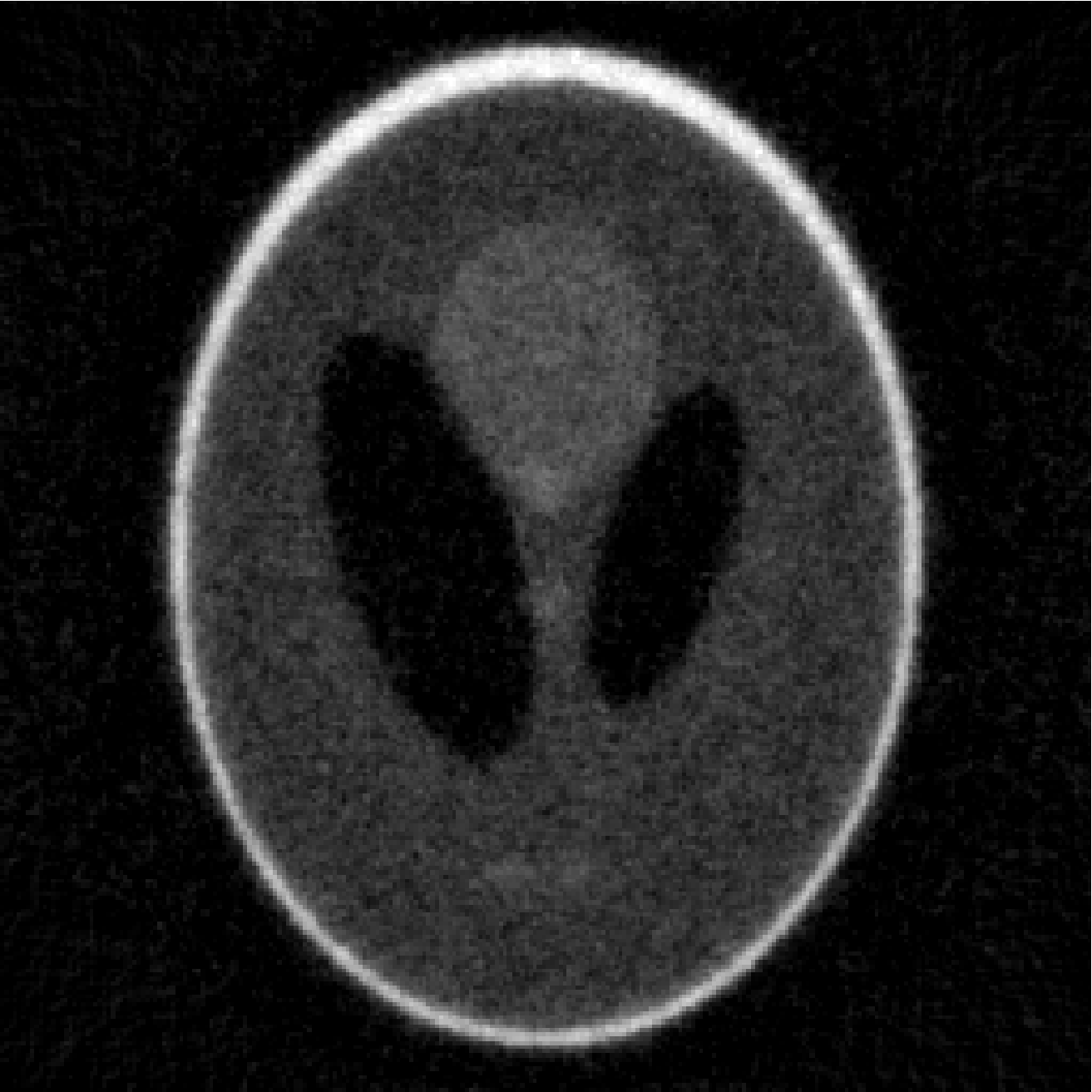}}  \\
   \subfloat[AB-GKB]{\label{Fig:Shepplogan ABGKBRNS}\includegraphics[width=.20\textwidth]{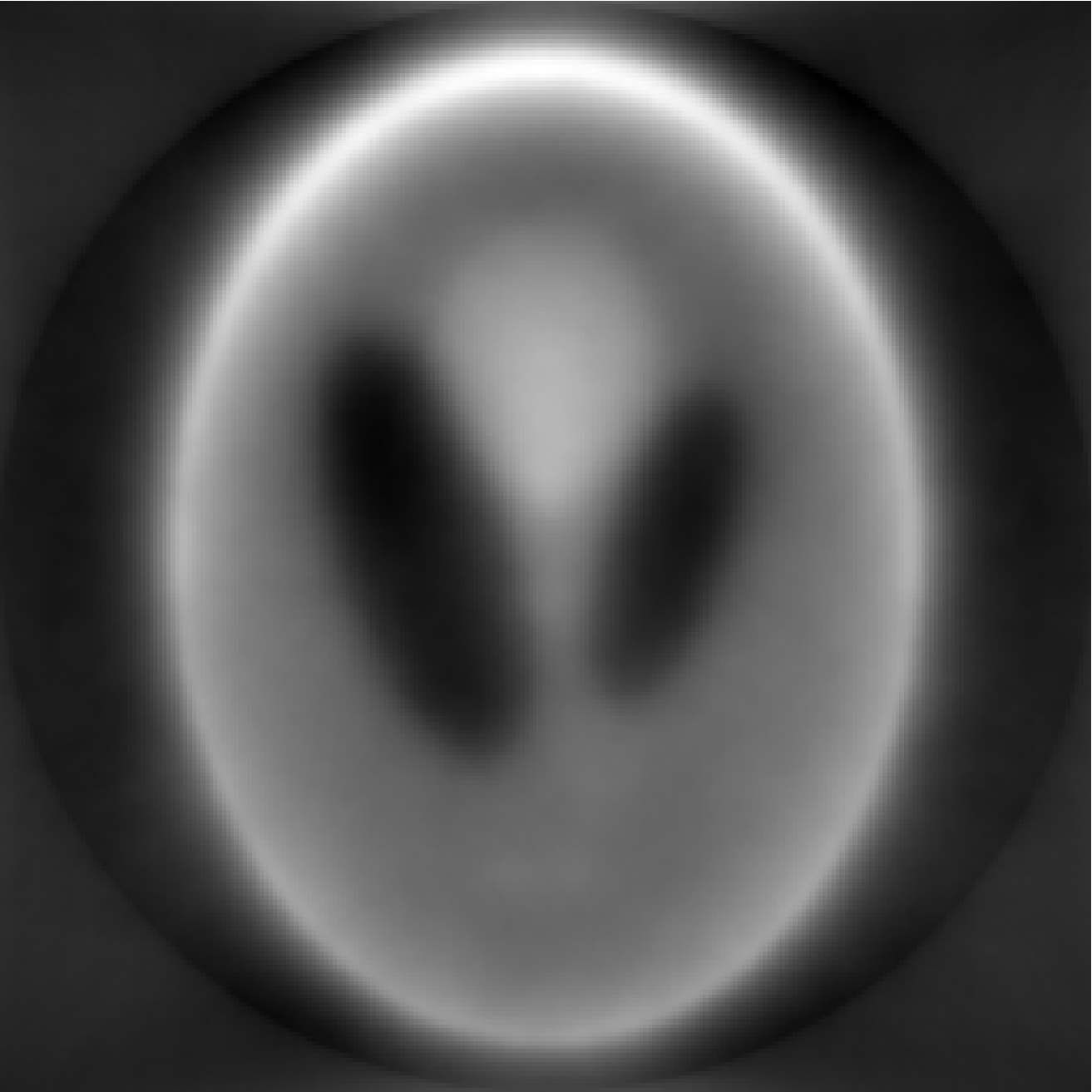}}   \ 
     \subfloat[BA-GKB]{\label{Fig:ShepploganShepplogan BAGKBRNS}\includegraphics[width=.20\textwidth]{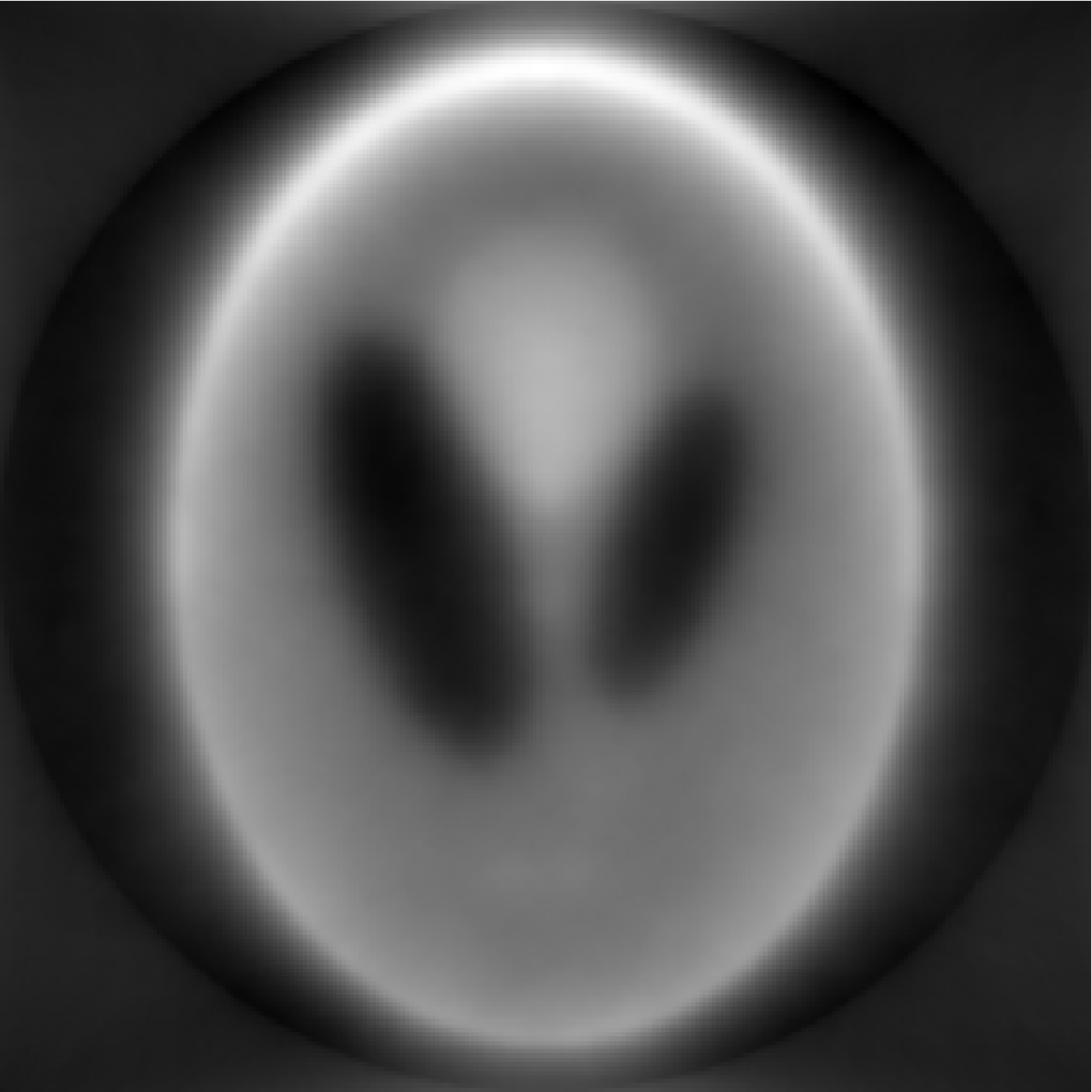}} \
  \subfloat[AB-GMRES]{\label{Fig:Shepplogan ABGNRESRNS}\includegraphics[width=.20\textwidth]{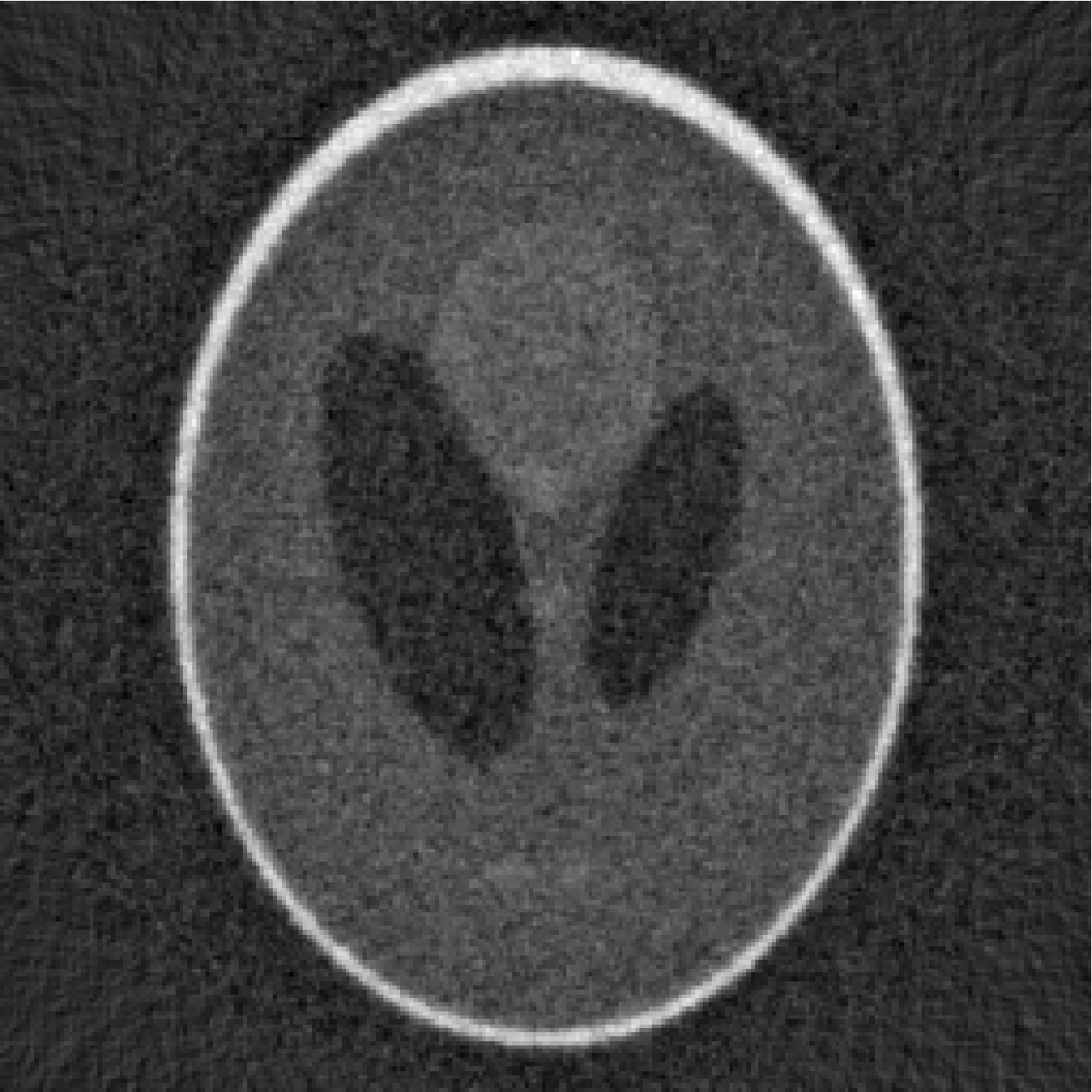}}  \
  \subfloat[BA-GMRES]{\label{Fig:Shepplogan BAGMRESRNS}\includegraphics[width=.20\textwidth]{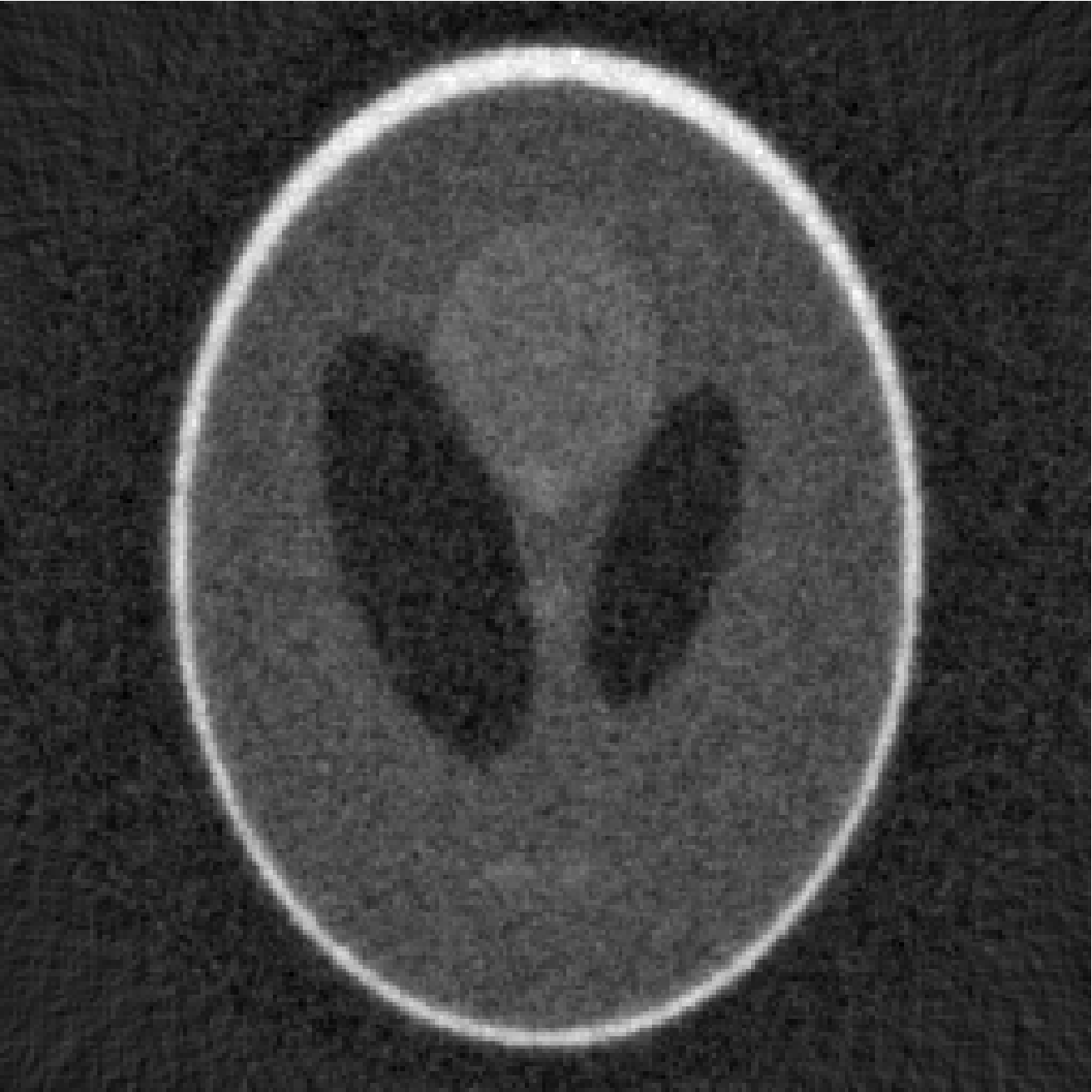}}  \
  \caption{Reconstructed images corresponding to \Cref{Fig:blurred Shepplogan}, obtained using DP and RNS stopping criteria in the first and second rows, respectively.}\label{Fig:Recons Shepplogan}  
  \end{figure}

In this problem, we use a large noise to evaluate the impact on the proposed algorithms. As shown in \Cref{Fig:Shepplogan RRE}, the AB- and BA-GKB methods generally achieve lower RRE than the AB- and BA-GMRES methods, except during the early iterations. Moreover, we see that GMRES approaches are more affected by semiconvergence, whereas the RRE obtained by GKB methods stabilizes after approximately $100$ iterations. The simulation also demonstrates that both GKB variants perform comparably. The corresponding reconstructions are presented in \Cref{Fig:Recons Shepplogan}.

The results reported in \Cref{Tab:Shepplogan Problem} show that the DP criterion effectively terminates all algorithms, yielding comparable RRE values. Moreover, the GMRES methods require a similar number of iterations with both DP and RNS to converge, indicating the effective choice for $\epsilon$. While the RNS criterion enables faster termination for the GKB methods, it leads to solutions with relatively higher RRE. Column $4$ of \Cref{Tab:Shepplogan Problem} further indicates that the GMRES methods are computationally less expensive than the GKB methods, which is consistent with the cost analysis presented in \Cref{sec:comp costs}. Moreover, for the underdetermined case, both BA-GKB and BA-GMRES are less expensive than their AB counterparts.

\begin{table}[htb]
 \footnotesize 
 \caption{\label{Tab:Shepplogan Problem} The corresponding \text{RRE} \eqref{eq:RE} computed at $k_{DP}$ \eqref{eq:DP} and $k_{RNS}$ \eqref{eq:RNS}, by all methods, applied to the problem with $\text{SNR}\approx 20$ in \Cref{ex:problem2}. The lowest RRE results are shown in boldface. Timings are recorded for $150$ iterations of all algorithms. The results obtained using the RNS for the
corresponding methods are provided in parentheses.}
 \begin{tabular}[t]{lccccccccc } 
 \toprule 
Method&Iteration&$\text{RRE}(\bfx)$&Time(s)  \\ \midrule
 AB-GKB&$21(6)$&$\boldsymbol{0.34}(0.60)$&$3.88$\\
 BA-GKB&$23(6)$&$\boldsymbol{0.34}(0.61)$&$4.26$\\  
 AB-GMRES&$6(7)$&$0.35(0.35)$&$2.85$\\ 
 BA-GMRES&$6(7)$&$0.36(\boldsymbol{0.34})$&$3.09$\\ 
\bottomrule
 \end{tabular}    
 \end{table}

Our numerical examples, which are not reported here, indicate that varying noise levels across all cases do not affect our conclusions, thereby confirming the robustness and reliability of the algorithms we proposed.

%%%%%%%%%%%%%%%%%%%%%%%%%%%%%%%%%%%%%%%%%%%%%%%%%%%%%%%%%%%%%%%%%%%%%%%%%%%%%%%%%%%

 \section{Conclusions}\label{sec:Conclusion}
In this work, we proposed two Krylov subspace solvers, namely AB-GKB and BA-GKB, designed for reconstructed large-scale X-ray computed tomography problems involving unmatched projector pairs. These methods leverage the flexibility of Krylov subspace techniques in conjunction with the SVD generated via GKB, thereby effectively mitigating the semiconvergence behavior commonly observed in classical Krylov methods. Achieving stable reconstructions, however, requires an appropriate stopping criterion to terminate iterations before noise amplification.

The theoretical analysis of computational complexity suggests that, while our proposed algorithms efficiently compute the SVD of a sparse lower bidiagonal matrix—rather than the SVD of a Hessenberg matrix as seen in GMRES methods—they incur a higher computational cost compared to AB-GMRES and BA-GMRES methods, which require fewer applications of forward and backward projectors. Nonetheless, our approach significantly reduces storage and computational demands by projecting the original problem onto a lower-dimensional subspace and solving the resulting problem via SVD.

Numerical experiments on computed tomography inverse problems demonstrate that the proposed methods yield stable reconstructions without compromising the accuracy achieved by GMRES methods. Furthermore, the numerical results highlight the effectiveness of the DP as a stopping rule, which outperforms the RNS criterion. Future work will focus on enhancing the efficiency and robustness of the proposed algorithms by incorporating advanced regularization techniques that further exploit the complementary strengths of the GKB framework.

\section*{Acknowledgments}
The author extends his appreciation to the Deanship of Scientific Research, Islamic University of Madinah, Saudi Arabia, for funding this research work.

\section*{Funding}
Not applicable.

\bibstyle{natbib}
\bibliography{sn-bibliography}

\end{document}